\documentclass[11pt]{article}
\usepackage{amsmath,amsthm,amssymb,amsfonts}
\usepackage{verbatim}
\usepackage{latexsym}
\usepackage{graphicx,psfrag,import}
\usepackage{fullpage}
\usepackage{framed}
\usepackage{verbatim}
\usepackage{color}
\usepackage{epsfig}
\usepackage{epstopdf}
\usepackage{hyperref}
\usepackage{geometry}
\usepackage{mathtools}
\usepackage{enumerate}
\usepackage{multicol}
\usepackage{a4wide}
\usepackage{booktabs}
\usepackage{enumitem}
\usepackage{lineno}
\usepackage{parcolumns}
\usepackage{thmtools}
\usepackage{xr}
\usepackage{epstopdf}
\usepackage{mathrsfs}
\usepackage{comment}
\usepackage{authblk}
\usepackage{setspace}
\usepackage[colorinlistoftodos]{todonotes}
\usepackage{physics}
\usepackage{stackrel}
\usepackage{stmaryrd}
\usepackage{subcaption}
\usepackage{float}

\definecolor{bittersweet}{rgb}{1.0, 0.44, 0.37}
\definecolor{brightturquoise}{rgb}{0.03, 0.91, 0.87}
\definecolor{burntsienna}{rgb}{0.91, 0.45, 0.32}
\definecolor{citrine}{rgb}{0.89, 0.82, 0.04}
\definecolor{paleplum}{rgb}{0.8, 0.6, 0.8}

\usepackage{adjustbox}

\def\id{\mathrm{id}}

\newtheorem{theorem}{Theorem}[section]
\newtheorem{corollary}[theorem]{Corollary}

\newtheorem{lemma}[theorem]{Lemma}
\newtheorem{remark}[theorem]{Remark}

\newtheorem{examplex}[theorem]{Example}
\newenvironment{example}
  {\pushQED{\qed}\examplex}
  {\popQED\endexamplex}

\theoremstyle{definition}
\newtheorem{definition}[theorem]{Definition}
\theoremstyle{definition}

\def\ds{\displaystyle}
\def\been{\begin{enumerate}}
\def\enen{\end{enumerate}}
\def\beit{\begin{itemize}}
\def\enit{\end{itemize}}

\def\Z{\mathbb{Z}}

\newcommand{\R}{\mathbb{R}}

\def\xrin{\xrightarrow{t \to \infty}}

\usepackage{tikz}
\usetikzlibrary{arrows,shapes.gates.logic.US,shapes.gates.logic.IEC,calc}
\usepackage{circuitikz}

 \tikzset{every node/.style={auto}}
 \tikzset{every state/.style={rectangle, minimum size=0pt, draw=none, font=\normalsize}}

\newcommand*\circled[1]{\tikz[baseline=(char.base)]{
\node[shape=circle,draw,inner sep=2pt] (char) {#1};}}
            
\ctikzset{logic ports origin=center}

\title{Chemical mass-action systems as  analog computers: \\  implementing arithmetic computations at specified speed}

\author{David F. Anderson\footnote{Department of Mathematics, University of Wisconsin-Madison, {\tt anderson@math.wisc.edu}} ~ and Badal Joshi\footnote{Department of Mathematics, California State University San Marcos, {\tt bjoshi@csusm.edu}.}}

\begin{document}

\maketitle

\begin{abstract} 
     Recent technological advances allow us to view chemical mass-action systems as analog computers. 
     In this context, the inputs to a computation are encoded as initial values of certain chemical species while the outputs are the limiting values of other chemical species. 
     In this paper, we design chemical systems that carry out the elementary arithmetic computations of: identification, inversion, $m$th roots (for $m \ge 2$), addition, multiplication, absolute difference, rectified subtraction over non-negative real numbers, and partial real inversion over real numbers.  We prove that these ``elementary modules'' have a speed of computation that is independent of the inputs to the computation.  
     Moreover, we prove that finite sequences of such elementary modules, running in parallel, can carry out composite arithmetic over real numbers, also at a rate that is independent of inputs.  
     Furthermore, we show that the speed of a composite computation is precisely the speed of the slowest elementary step. 
     Specifically, the scale of the composite computation, i.e. the number of elementary steps involved in the composite, does not affect the overall asymptotic speed -- a feature of the parallel computing nature of our algorithm.
     Our proofs require the careful mathematical analysis of certain non-autonomous systems, and we believe this analysis will be useful in different areas of applied mathematics, dynamical systems, and the theory of computation.  
     We close with a discussion on future research directions, including numerous important open theoretical questions pertaining to the field of computation with reaction networks.
     
     {\bf Keywords:} analog computation, computing with chemistry, polynomial dynamical systems, input-independent speed of computation
     
     {\bf MSC:} 37N25, 68N01, 92B05
\end{abstract}

\section{Introduction}

An increasingly active area of scientific research is the study of chemical reaction systems as computational machines \cite{anderson2021reaction, buisman2009computing, cappelletti2020stochastic, chen2023rate, chen2014deterministic, qian2011neural, QSW2011, qian2011simple, SSW2010}.  
The broad goal of this nascent field is to develop systems that can operate in the niche of a (wet) cellular environment, rather than to directly compete with modern digital computers.
Methods for the physical implementation of these systems include DNA-strand displacement,  which was developed largely by  Lulu Qian, David Soloveichik, and Erik Winfree \cite{QSW2011,SSW2010}.  While this research area is still in its relative infancy, it is critical that the theoretical (mathematical) foundation of these methods and systems keeps pace with the physical implementations so as to (i) keep the algorithmic aspects of this field on a firm mathematical footing from the beginning, and (ii)  provide a reliable (proof based) method for the development of the appropriate reaction network-based software design.

Under the assumptions that the chemical system is well-stirred and at a constant temperature, deterministic mass-action chemical systems are mathematically modeled as autonomous systems of ordinary differential equations \cite{feinberg2019foundations}.  
To utilize such a system as a computational machine, we associate the inputs and outputs of the desired computation with certain state variables of the dynamical system. Most often, the desired output is encoded as a limiting value of a subset of the state variables  that converge asymptotically to this limiting value. Since by necessity any recording of the system is at a finite time it is important to bound the error, which in this case is the deviation from the asymptotic value. The rate of this  convergence, or alternatively the rate at which the error decays, determines the speed of the computation.   The time chosen for the recording will in general be determined by the required level of accuracy, but should not be so large as to render the method impractical.

Because mass-action systems are modeled via systems of differential equations, they can naturally be viewed as analog computers. 
While digital computers now dominate the computational landscape in our society, it is interesting to note that analog computers have certain advantages over digital computers.  For example, consider the problem of numerically solving the initial value problem
\begin{equation}\label{eq:6785785}
    \dot x(t) = 1 - x(t), \quad x(0) = x_0 \in \R_{>0},
\end{equation}
If you were using a digital computer, you would naturally utilize some sort of discretization method (such as Euler's method or a Runge-Kutta method) to approximate the solution to some desired accuracy. No matter the method, however, the key point is that one \textit{must} discretize,  as digital computers can only handle discrete computations.  An incredible amount of scientific effort has gone into discretization methods for the numerical solution of ODEs, PDEs, and SDEs via digital computers over the past many decades, with such work still dominating portions of applied mathematics.  

On the other hand, the system \eqref{eq:6785785} above is quite naturally ``solved'' by simply considering the chemical reaction network modeled via
\[
\emptyset \rightleftarrows X,
\]
with an initial concentration of $x_0$, so long as the rate constants are chosen to be equal to one (see section \ref{sec:mass_action_ODEs} for the details of mass-action models).  In short: many problems that are difficult to solve for digital computers are quite simple to solve for a chemical system.

However, there is a natural converse to the above observation.  While chemical systems will naturally be able to solve systems that are  described via differential equations (see section \ref{sec:mass_action_ODEs} for the exact class of differential equations they solve), it is more difficult to use them to compute things that are naturally discrete, such as basic arithmetic (addition, multiplication, finding roots, division, etc.).  Hence, while it took much thought by many people to develop algorithmic methods for the solution of continuous problems with digital computers, it is natural that we will require a fair amount of work to develop efficient methods for the use of chemical systems to solve these discrete problems.

We therefore focus here on algorithmic  aspects of chemistry-based analog computing of arithmetic operations, as these operations are a cornerstone of discrete computation.  Specifically, we focus on the development of elementary ``chemical modules'' that  carry out  basic arithmetic operations (such as inversion, roots, addition, multiplication, etc.) and we prove that the developed modules have a speed of computation that is independent of the inputs.\footnote{More specifically, we prove that the speed of computation for the elementary modules have a lower bound that is independent of inputs.}  Moreover, we prove that these elementary modules can be run in parallel so as to carry out any desired arithmetic computation, also at speed that is independent of input.

There is existing work on algorithmic aspects of analog computation of arithmetic operations, even from the specific vantage point of reaction network implementation. 
Our work is closest in approach to the one taken by \cite{buisman2009computing}. 
However, while the authors of \cite{buisman2009computing} considered algorithms for arithmetic computations using reaction networks, they did not include \textit{speed of computation} in their consideration. 
Our stated purpose is to remedy this. 
We frame a general arithmetic computation as a finite sequence of elementary computations that are running in parallel. Since every element of the sequence is being computed in parallel, we need to consider time-evolving inputs to the elementary modules, where the output of a previous step is providing input into the next step even without having converged to its asymptotic value. 
We show input-independence of speed of any  overall computation that is composed from multiple, but finitely many, elementary arithmetic operations.

Throughout this work, we assume that the rate constants of our models (which are the key system parameters, see section \ref{sec:mass_action_ODEs}) are all exactly equal to 1.  Our main results of this paper, Theorem \ref{thm:composite} and Corollary \ref{cor:speed_special}, then concludes that any arithmetic computation carried out in the manner specified herein, even those consisting of a finite number of composite calculations being carried out in parallel, has a rate of convergence of at least 1 (loosely, this means the deviation of the model from   the desired output decreases like $e^{-t}$). 
While it is important to our work that all the rate constants be equal, the choice of setting them equal to 1 is arbitrary.  We discuss in Remark \ref{remark:scaling} how setting them all to $\sigma>0$ is simply a time-change of our model and leads to a rate of convergence of $\sigma$.   We discuss the possibility of generalizing to the situation of non-equal rate constants in section \ref{sec:discussion}.

There are several bodies of research, outside the world of chemical reactions, on designing efficient hardware and software for analog as well as analog-digital hybrid computing \cite{tsividis2018not,maclennan2014promise,achour2016configuration,daniel2013synthetic}. 
While we develop our results for a chemistry-based computer, the results apply equally well to any computational system that computes via nonlinear polynomial differential equations. We therefore hope that our theoretical results will find use in the broader sphere of computation.

The remainder of the paper is organized as follows.  In the next subsection, section \ref{sec:comparing}, we point out that there are two natural ways to use chemical systems as computational machines and compare/contrast their strengths and weaknesses.  In section \ref{sec:BasicProblem}, we (i) formally introduce the mathematical models that will serve as our algorithms for analog computation, which are deterministic chemical reaction networks with mass-action kinetics, (ii) discuss issues related to these models, and (iii) demonstrate, via example, the central purpose of this paper: the development of reaction network modules that can carry out arithmetic at speeds that are independent of inputs.  In section \ref{sec:math}, we provide the main mathematical analysis of this paper.  In particular, we provide a rather detailed analysis of certain non-autonomous dynamical systems. The work of section \ref{sec:math} is general enough that we hope and expect it to have applications in other areas of applied mathematics. In section \ref{sec:elem_comp}, we provide the chemical reaction network modules that carry out the basic elementary operations discussed above (identification, addition, multiplication, inversion, etc.) and use the mathematical results of section \ref{sec:math} to prove that they have computational speeds that are independent of inputs.  In section \ref{sec:combining}, we demonstrate how any calculation carried out via a combination of elementary operations using our constructions also has a computational speed that is independent of inputs.  In section \ref{sec:reals}, we demonstrate how to expand the constructions to handle calculations involving all real numbers (as opposed to non-negative numbers).  Finally, in section \ref{sec:discussion} we provide a discussion, including many directions for future research.

\subsection{Arithmetic with reaction networks: comparing existing schemata}
\label{sec:comparing}

There are at least two  schemata that have been used for carrying out arithmetic with reaction networks.  We refer to the first schema as {\em conservation law-based}, with examples found in \cite{chen2023rate, song2016analog}.  See also \cite{salehi2018computing}, which uses a {\em fractional  representation} to represent real numbers in the interval $[0,1]$.  General functions can then be implemented via the methods in \cite{salehi2018computing} by  first approximating the desired function by a Taylor series and then implementing the resulting polynomial function with the relevant reaction network.  The second schema, which we refer to as  {\em positive steady state-based} is the focus of the current paper (see also \cite{buisman2009computing}).

To understand one way that addition could be performed in the {\em conservation law schema}, we consider a simple analogy. 
Suppose we have two `input' buckets and one `output' bucket. 
The contents of the buckets may be discrete (e.g., marbles) or continuous (e.g., water). 
Pouring the contents of the two input buckets into an initially empty output bucket will result in the sum of the two inputs as the eventual volume of the output bucket. 
A reaction network based on this simple idea is 
\begin{align} \label{net:add0}
    A \to X, \quad B \to X. 
\end{align}
Here $A$ and $B$ are input species while $X$ is the output species. Clearly, if the initial population of $X$, represented by $x$, is $0$, we will have that, asymptotically, $x$ converges to the sum of the initial populations of the input species: 
\begin{align} \label{eq:add0_de}
    \lim_{t \to \infty} x(t) = a(0) + b(0). 
\end{align}
This is true regardless of the specific choice of rate constants for the chemical reactions (see section \ref{sec:mass_action_ODEs}).
 In fact, it is true independent of whether the underlying model is discrete (integer counts) or continuous (concentrations), deterministic or stochastic. 
Consequently, the conservation law based schema has the potential to be used in both digital as well as analog computing.

Despite the simplicity and potential use in both digital- and analog-computing, the conservation law based schema does have significant drawbacks which we now discuss. 
Perhaps the most significant one is the `loss of memory' suffered by the input species. 
At the end of the addition computation, all memory of what the inputs contained is lost forever. 
In fact, if addition is one of several elementary steps in a longer sequence of computations, then all inputs and all intermediate steps are forgotten. Only the final output is retained. 
This may become an issue if the same input is required multiple times. 
There are methods for remedying this. One idea is to introduce back-up species which re-record the lost input. For instance, the reaction network for addition \eqref{net:add0} may be replaced by the following reaction network  
\begin{align} \label{net:add1}
    A \to \widetilde{A} + X, \quad B \to \widetilde{B} + X. 
\end{align}
Here the species $\widetilde{A}$ and $\widetilde{B}$ back up the input species $A$ and $B$, because in the long run
\begin{align} \label{eq:add1_de}
    \lim_{t \to \infty} \widetilde{a}(t) = a(0),  \quad \lim_{t \to \infty} \widetilde{b}(t) = b(0), \quad \mbox{ and } \lim_{t \to \infty} x(t) = a(0) + b(0),
\end{align}
so long as we began with $\widetilde{a}(0) = \widetilde{b}(0) = 0$.

The {\em positive steady state based schema} utilized in the present paper does not suffer from loss of inputs, and so separate back up species are not required at every step. 
For every elementary operation, we explicitly require that all inputs are preserved while the output is computed and stored as the concentration of a dedicated new species. 

To illustrate the `positive steady state based' schema utilized in this paper let us once again consider the operation of addition, now accomplished with the reaction network
\begin{align} \label{eq:add12}
    A \to A + X, \quad B \to B + X, \quad X \to 0.
\end{align}
Assuming rate constants that are equal to one, the mass action differential equation for these reactions is (see section \ref{sec:mass_action_ODEs})
\begin{align*}
    \dot x(t) = a(t) + b(t) - x(t). 
\end{align*}
We first consider the case where the reactions in \eqref{eq:add12} stand in isolation. 
Note that the concentrations of $A$ and $B$ are not changed by any of the reactions in \eqref{eq:add12}. 
In particular, this means that the input species $A$ and $B$ retain their initial values, $a(t) = a(0)$ and $b(t) = b(0)$ for all time $t \ge 0$. 
Meanwhile, the concentration of the output species $X$ asymptotically converges to $a(0) + b(0)$. 
Now imagine that the reactions in \eqref{eq:add12} are only some reactions in a much larger network. 
The inputs $A$ and $B$ may be outputs of some other subnetwork and may as yet be in the process of converging to their asymptotic values. 
Then, assuming the limits on the right exist, we have
\begin{align*}
    \lim_{t \to \infty} x(t) = \lim_{t \to \infty} \left( a(t) + b(t) \right).
\end{align*}
This example illustrates how an output from one computation may be fed into another computation as well as how inputs may be simultaneously used in several different computations since they are not degraded.

\section{Introducing the basic problem}
\label{sec:BasicProblem}

In this paper, we consider the task of carrying out the following common arithmetic functions via deterministic mass-action chemical systems:  addition, multiplication, inversion, division, $m$th root (for integers $m \ge 2$), absolute difference of two numbers, the maximum of two numbers, and subtraction.  Our particular goal, which we will describe more fully below,  is to provide mass-action chemical reaction network modules for each such operation whose ``speed of computation'' is independent of the inputs (and hence outputs) to the problem \textit{and} which can be combined to run in parallel so that any sequence of arithmetic computations (e.g., the square root of the inverse of the sum of two numbers: $\sqrt{1/(a+b)}$) will have a speed of computation that is also independent of inputs.

There are a number of questions/issues we wish to highlight from the start and are the focus of the present section. 
\begin{enumerate}
    \item What is the class of differential equation models we use in this paper? Said differently, what class of differential equations arise from deterministically modeled mass-action chemical systems?

    \item Chemical concentrations are naturally non-negative. 
 Hence, it is natural to consider problems such as the addition or division of two positive reals since the inputs and outputs are positive reals.  However, how should one handle the situation in which one or more of the inputs are negative?  Similarly, how should one handle the operation of minus, whose output can be negative even when the inputs are both positive?

 \item Is it possible to ensure that the speed of the computation does not depend on the inputs and/or outputs to the problem?  It is the recognition of this problem and its resolution in the affirmative that we consider to be the main contribution of this paper.
 
\end{enumerate}

We handle each of the above questions in turn before moving  to our general theory and constructions.
Before proceeding, however, we first note that efficiency, accuracy, and speed are essential, pragmatic considerations when implementing computers. 
We need some assumptions, however, regarding the interface with a chemistry-based ``CPU.'' 
Specifically, we assume that we have (i) the ability to initiate chemical species with arbitrary precision, and (ii) the ability to measure chemical species with arbitrary precision.

\subsection{Mass-action chemical systems and their associated ODEs}
\label{sec:mass_action_ODEs}

Consider an ordinary differential equation of the form
\[
\dot x(t) = f(x(t)),
\]
with $x(t) \in \R^d$ and $f:\R^d \to \R^d$.  When does such an equation arise from a mass-action model of a chemical system?

Before answering, we first introduce mass-action models.  Consider a finite number of chemical ``species,'' which we will denote by $\{X_1,\dots, X_d\}$, undergoing transformations via a finite number of possible reaction types.  Reactions can be understood intuitively: they require certain numbers of the different species as inputs, and they output a different number of the various species.  They are most naturally visualized in graphical format.  For example, consider a system with three species, $\{X_1,X_2,X_3\}$, and the following three reaction types
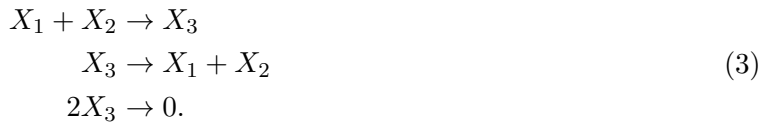
\begin{align}
\begin{split}\label{eq:7896}
    X_1 + X_2 &\to X_3\\
    X_3 &\to X_1 + X_2\\
    2X_3 &\to 0.
    \end{split}
\end{align}
In the first reaction, one molecule of $X_1$ and one molecule of $X_2$ are converted to one molecule of $X_3$.  The second reaction reverses the first.  The third reaction consists of the destruction of two  $X_3$ molecules (perhaps they merge to become a fourth type of species that we are not modeling).  In general, a chemical system will have a finite number of such reactions of the form
\begin{equation}
    \label{eq:0986545678}
\nu\cdot X \to \nu' \cdot X,
\end{equation}
where $\nu,\nu'\in \Z^d_{\ge 0}$ and by $\nu\cdot X$ we mean the dot product:  $\nu\cdot X = \nu_1X_1 + \cdots + \nu_d X_d$.  Note that if we were explicitly counting the numbers of molecules of each species, then one instance of the  reaction \eqref{eq:0986545678} would change the system by addition of the vector $\nu'-\nu\in \Z^d$.

 A differential equation model for the chemical concentrations arises by associating to each  reaction a rate function.  For example, if we denote the above  reaction \eqref{eq:0986545678} by $\nu\to \nu'$ and denote the rate of that reaction when the system is in state $x\in \R^d_{\ge 0}$ by $\lambda_{\nu\to \nu'}(x)$, then the differential equation governing the dynamics is
\begin{equation}
\label{eq:mass-actionODE}
\dot x(t) = \sum_{\nu\to \nu'} \lambda_{\nu\to \nu'}(x(t)) (\nu'-\nu),
\end{equation}
where the sum is over all reaction types (enumerated by $\nu\to \nu'$).
The system is said to be governed by \textit{mass-action} kinetics when the functions $\{\lambda_{\nu\to\nu'}\}$ have the following form,
\[
\lambda_{\nu\to \nu'}(x) = c_{\nu\to \nu'} \prod_{i=1}^d x_i^{\nu_i},
\]
where each $c_{\nu\to \nu'}$ is a positive constant and where we interpret $0^0$ to be $1$.  Hence, for mass-action kinetics the rate of a reaction is a monomial in which the exponent for the $i$th species is the coefficient of $X_i$ in $\nu\cdot X$.   For example, assuming mass-action kinetics the differential equations governing the system \eqref{eq:7896} are
\begin{align*}
    \dot x(t) = c_{1} x_1x_2 
    \left( \begin{array}{r}
    -1\\
    -1\\
    1
    \end{array}
    \right) + c_{2} x_3
    \left( \begin{array}{r}
    1\\
    1\\
    -1
    \end{array}
    \right) + c_{3} x_3^2
    \left( \begin{array}{r}
    0\\
    0\\
    -2
    \end{array}
    \right),
\end{align*}
where  $c_{1}, c_{2},$ and $c_{3}$ are some positive constants, and where we enumerated the constants sequentially (as opposed to the reaction $\nu\to \nu')$).  It is common to incorporate the rate constants into the network graph, such as \eqref{eq:7896}, by placing them next to the associated arrow.

We can now answer the first question posed above.  If $f: \R^d \to \R$ satisfies the two properties below, then there is a mass-action chemical system whose associated ODE \eqref{eq:mass-actionODE} is given by $\dot x(t) = f(x(t))$.

\begin{enumerate}
\item The function $f$ is a polynomial.  That is, for each $i \in \{1,\dots,d\}$, 
\begin{equation}\label{eq:mass-actionsystem}
f_i(x) = \sum_{a_1 = 0}^N\cdots \sum_{a_d=0}^N c_{i,a_1,\dots,a_d}x_1^{a_{1}}\cdots x_d^{a_{d}},
\end{equation}
where $N<\infty$ is a maximum exponent in the model, and the $c_{i,a_1,\dots,a_d}$ are constants (which can be negative or zero). 

\item If $c_{i,a_1,\dots, a_d}<0$ in  \eqref{eq:mass-actionsystem}, then $a_i\ge 1$ for the associated monomial.  In plain English, $x_i$ must appear in any monomial of $f_i(x)$ that has a negative coefficient.
\end{enumerate}
Note that the final condition guarantees that the non-negative orthant $\R^d_{\ge 0}$ is forward invariant under the dynamics.  That is, so long as $x(0) \in \R^d_{\ge 0}$, then $x(t) \in \R^d_{\ge 0}$ for all $t \ge 0$.

It should be clear that any mass-action system satisfies the two conditions above.  To see the converse, that is to see that for any dynamical system $\dot x(t) = f(x(t))$ satisfying the conditions above there is a mass-action chemical system with $f$ as the governing dynamics, see  \cite{hars1981inverse} or, for a more modern treatment, \cite{feinberg2019foundations}.   

\subsection{The issue of negativity and the dual rail representation} \label{sec:negativity}

Chemical concentrations are non-negative, but we wish to be able to represent negative values within our framework.  To do so, we represent real numbers as the difference between two given chemical species.  This is sometimes referred to as the \textit{dual rail representation} of real numbers.  The introduction we give to this topic in this section will be brief.  See, for example, \cite{chen2023rate} for a more thorough discussion. 

One of the early works to utilize the dual rail representation is \cite{chen2014deterministic}, where, for example, they labeled $Y^P$ and $Y^C$ the amounts that a particular species, $Y$,  had been ``P''roduced and ``C''onsumed, respectively.  Then, the total amount of $Y$ was $Y^P-Y^C$.    In this example, if, for example, $Y^P = 10$ and $Y^C = 8$, then the total amount of $Y$ is clearly $2$.  However, we could also represent the value 2 by having $Y^P= 2+c$ and $Y^C = c$ for any $c \ge 0$. Thus, \textit{how} we represent a given real number is a bit ambiguous, a priori.

 For our purposes, we wish to do away with such ambiguity.   Hence, we will split a number into its positive and negative parts.   That is,  we define the {\it dual rail map} $\varphi: \R \to \R_{\ge 0}^2$ via
\[
a \mapsto (a_p, a_n) = 
\begin{cases}
(a,0),  & \mbox{if } a > 0, \\ 
(0,-a), & \mbox{if } a \le 0. 
\end{cases}
\]
Here $a_p$ and $a_n$ are referred to as the {\it positive} and {\it negative parts}, respectively, of $a$. 
The map $\varphi$ is injective and 
\[
\Im(\varphi) = \{(a_p, a_n) \in \R^2_{\ge 0} ~|~ a_p a_n = 0\} = \{(a_p, a_n) \in \R^2_{\ge 0} ~|~ a_p=0 \text{ or } a_n = 0\}. 
\]
Define the map $\psi:  \R_{\ge 0}^2 \to \R$ via $(a_p, a_n) \mapsto a = a_p - a_n$. 
The map resulting from restricting the domain of $\psi$ to $\Im(\varphi)$ is the inverse of $\varphi$.  
We will refer to $(a_p,a_n)$ as the canonical representation of $a$ when $a_pa_n=0$, and write $a = (a_p, a_n)$, for brevity.

\subsection{Input independent speed of computation} \label{sec:input_independent}

Suppose we want to compute the reciprocal of a positive real $a \in \R_{> 0}$ as the long-term solution of a mass-action differential equation. 
A natural attempt would be to consider the reaction network
\begin{align} \label{net:inv_naive}
0 \to X, \quad A+X \to A
\end{align}
whose mass-action differential equation, assuming all reaction rate constants are equal to $1$ and that the concentration of species $A$ is $a$, is 
\begin{align} \label{ode:inv_naive}
    \dot x = 1 - a x. 
\end{align}
The differential equation has a unique solution for all time given by 
\[
    x(t) = \frac{1}{a} + \left(x(0) - \frac{1}{a}\right) e^{-at}.
\]
Clearly, $x(t) \xrin 1/a$. 
Moreover, since $a$ is the coefficient of $t$ in the exponential decay $e^{-at}$, the rate of convergence to the long term solution is $a$. 
The rate of convergence is directly related to the speed of computing the reciprocal. 
Indeed, suppose that $T_n$ is the smallest time required to compute $1/a$ accurately up to $n-1$ digits. Then we have that 
\begin{align*}
   \abs{x(T_n) - 1/a}   =  \abs{x(0) - 1/a}e^{-a T_n} = 10^{-n}, 
\end{align*}
which after applying logarithms and rearranging gives us 
\begin{align*}
T_n = \frac1a \ln\left( \abs{x(0) - 1/a} \right) + \frac{n \ln(10)}{a}. 
\end{align*}
In particular, we have that 
\begin{align} \label{eq:limittime0}
\lim_{n \to \infty} \frac{T_n}{n \ln(10)} = \frac1a, 
\end{align}
so that the time required to compute $n$ decimal digits of $1/a$ scales not only with $n$ but with the output $1/a$. 
Alternatively, the reciprocal 
\begin{align} \label{eq:limitspeed0}
\lim_{n \to \infty} \frac{n \ln(10)}{T_n} = a, 
\end{align}
gives the speed of computation. We will say that the variable $x$  \textit{computes the reciprocal of $a$ at speed $a$}. 

We will now show that we can remove this dependence of speed of computation of the reciprocal on the input (and output) to the system. Then, over the rest of the paper, we will show that this input-independence can be achieved for an arbitrary computation. 

Consider the following reaction network  
\begin{align} \label{net:inv}
X \to 2X, \quad A+2X \to A + X
\end{align}
whose mass-action differential equation, assuming rate constants are equal to 1 and that the concentration of species $A$ is $a$, is
\begin{align} \label{ode:inv}
\dot x = x(1 -a x ). 
\end{align}
For $x(0) = x_0$, the unique solution for all time $t \ge 0$ is given by 
\begin{align*}
x(t) = \frac{x_0/a}{(1/a-x_0)e^{-t} + x_0}, 
\end{align*}
and so the distance from the attracting steady state $1/a$ is 
\begin{align*}
\abs{x(t)-1/a}  
= \left|  \frac{x_0/a}{(1/a-x_0)e^{-t} + x_0} - \frac1a\right|
=  \frac{1}{a} \left(\frac{\abs{x_0-1/a}}{x_0 e^t - (x_0-1/a)}\right). 
\end{align*}
This converges to zero like $e^{-t}$, and after a brief calculation, it is straightforward to see that 
\begin{equation} \label{eq:limittime1}
\lim_{n \to \infty} \frac{n \ln (10)}{T_n} = 1, 
\end{equation}
where, again, $T_n$ is the smallest time needed to compute $1/a$ to $n-1$ decimal places.  That is, $x$ computes the reciprocal of $a$ at speed 1.  A comparison of equations \eqref{eq:limittime0} and \eqref{eq:limittime1} highlights the input independence of the reciprocal computation in the latter construction.

\section{Rate of convergence for certain non-autonomous systems}
\label{sec:math}

In this section, which holds much of the mathematics, we provide a rather detailed analysis of certain non-autonomous dynamical systems.  These systems will have  a structure like
\[
\dot x(t) = f(g_1(t),g_2(t),x(t)),
\]
where $g_1$ and $g_2$ are best viewed as forcing functions that are themselves the output from another such system.  We will prove that for the particular systems we consider the following holds:  global exponential convergence of the forcing functions, such as $g_1$ and $g_2$, imply that the system solution also enjoys global exponential convergence.  That is, that property is  ``fed forward'' from the forcing functions
 to the solution.  We will also demonstrate that, in most cases, the rate of convergence for the system solution is at least as large as the minimum rate of convergence of the forcing functions.

To make these notions precise we provide here a definition and a lemma.  
Roughly speaking, ``rate of convergence'' of a function $g$ to a value $g^*$ is $\rho$ if $g(t) - g^* = f(t) e^{-\rho t}$, with $f(t)$ growing at most sub-exponentially, as $t\to \infty$. 
 Below, the reader may notice some connection with the commonly used ``big O'' notation.  However, we avoid this notation since it is necessary to be precise and make fine distinctions such as between ``rate that is at least $\rho$'' and ``rate that is greater than $\rho$''.

\begin{definition} \label{def:ratespeed}
Let $g: \R_{\ge 0} \to \R$ be a real-valued function that converges to a real number $g^* \in \R$. 
The \textit{rate of convergence} of $g$ to $g^*$ is defined to be 
\begin{align} \label{eq:ratespeed2}
\rho = - \limsup_{t \to \infty} \frac{\ln \abs{g(t) - g^*}}{t} 
\end{align}
whenever $\rho \in (0, \infty]$.
Alternatively, in the context of computation, we say that $g(t)$ computes $g^*$ at \textit{speed} $\rho$. 
\end{definition}

\begin{lemma} \label{lem:ratespeed}
Let $g: \R_{\ge 0} \to \R$ be a real-valued function that converges to a real number $g^* \in \R$.  Suppose there exist $c > 0$,  $\rho > 0$, and $m \ge 0$ such that for all $t \ge 0$
\[
\abs{g(t)-g^*} \le c t^m e^{-\rho t}, 
\]
then $g$ converges to $g^*$ at a rate that is at least $\rho$. 

Conversely, if $g$ converges to $g^*$ at a rate greater than $\rho$ then there exists a $c>0$ such that for all $t \ge 0$
\[
\abs{g(t)-g^*} \le c e^{-\rho t}.
\]
\end{lemma}
\begin{proof}
To see the first conclusion, take logarithms, divide by $t$, apply limit superiors to both sides and apply Definition \ref{def:ratespeed}.

For the converse, suppose that $g$ converges to $g^*$ at rate $\rho_g > \rho$, and let $\varepsilon = \rho_g - \rho$. 
Then there is a $T>0$ such that for any $t \ge T$, we have 
\begin{equation}
    \frac{\ln|g(t)-g^*|}{t} < -\rho_g + \varepsilon = -\rho. 
\end{equation}
Hence, for $t\ge T$, we have $|g(t)-g^*| < e^{-\rho t}$.  Defining
\[
c = \left(1 + \max_{t \in [0,T]}\{|g(t) - g^*|\} \right) e^{\rho T}, 
\]
we have $|g(t)-g^*| \le ce^{-\rho t}$ for all $t\ge 0$. 
\end{proof}

\subsection{Simple convergence results}

We begin with three lemmas.
\begin{lemma} \label{lem:sumproduct}
For each of $i \in \{1,2\}$, let $g_i:\R_{\ge 0} \to \R$ be a real-valued function that converges to $g_i^* \in \R$ at a rate that is at least $\rho_{g_i}$.
Then the following hold. 
\been
\item For $i \in \{1,2\}$, $\alpha \in \R \setminus \{0\}$, the scalar product $\alpha g_i$ converges to $\alpha g_i^*$ at rate that is at least $\rho_{g_i}$. 
\item The sum $g_1(t) + g_2(t)$ converges to $g_1^* + g_2^*$ at rate that is at least
\begin{align}
\min\{\rho_{g_1}, \rho_{g_2}\}. 
\end{align}

\item The product $g_1(t)g_2(t)$ converges to $g_1^*g_2^*$ at rate that is at least
\begin{align}
\begin{cases}
\min\{\rho_{g_1}, \rho_{g_2}\}, &\mbox{ if } g_1^* \ne 0,  g_2^* \ne 0, \\
\rho_{g_1}, &\mbox{ if } g_1^* = 0,  g_2^* \ne 0, \\
\rho_{g_2}, &\mbox{ if } g_1^*  \ne 0,  g_2^* = 0, \\
\rho_{g_1} + \rho_{g_2}, &\mbox{ if } g_1^* = 0,  g_2^* = 0. 
\end{cases}
\end{align}
\enen
\end{lemma}
\begin{proof}
We first prove the conclusions assuming that the rate of convergence of the functions $g_i$ is greater than $\rho_{g_i}$.
In that case, Lemma \ref{lem:ratespeed} implies that there exist positive constants $c_1$ and $c_2$ such that 
\[
\abs{g_i(t) - g_i^*} \le c_i e^{-\rho_{g_i} t}, \quad i \in \{1,2\}. 
\]
For the scalar product, we get the desired result from  
\[
\abs{\alpha g_i(t) - \alpha g_i^*} \le \abs{\alpha} c_i e^{-\rho_{g_i} t}, \quad i \in \{1,2\}, 
\]
and then applying the first part of Lemma \ref{lem:ratespeed}. 

Applying the triangle inequality to the sum of $g_1$ and $g_2$, we have 
\begin{align*}
\abs{g_1(t) +  g_2(t) - ( g_1^* +  g_2^*)} &\le  \abs{g_1(t) - g_1^*}  +  \abs{g_2(t) - g_2^*} \\
&\le c_1 e^{-\rho_{g_1} t} + c_2 e^{-\rho_{g_2} t}, 
\end{align*}
and so the result again follows from Lemma \ref{lem:ratespeed}. 

Turning to the product, we note first that for $i \in \{1,2\}$, 
\[
\abs{g_i(t)} - \abs{g_i^*} \le \abs{\abs{g_i(t)} - \abs{g_i^*}} \le \abs{g_i(t) - g_i^*} \le c_i e^{-\rho_{g_i} t}, 
\]
and so 
\[
\abs{g_i(t)} \le \abs{g_i^*} + c_i e^{-\rho_{g_i} t}. 
\]
Now for the product $g_1(t) g_2(t)$ consider,        
    \begin{align*}
\abs{g_1(t) g_2(t) - g_1^* g_2^*} &\le \abs{g_1(t) g_2(t) - g_1^* g_2(t)} + \abs{g_1^* g_2(t)- g_1^* g_2^*} \\
&= \abs{g_2(t)}\abs{g_1(t)  - g_1^* } + \abs{g_1^*} \abs{g_2(t)- g_2^*} \\
&\le \left(\abs{g_2^*} + c_2 e^{-\rho_{g_2} t}\right) c_1 e^{-\rho_{g_1} t} + \abs{g_1^*} c_2 e^{-\rho_{g_2} t}, 
    \end{align*}
and therefore the result follows.

To complete the proof, we now assume that one or both of the functions $g_i$ converge to their limiting values exactly at rates $\rho_{g_i}$. Then for an arbitrary $\varepsilon > 0$, the functions $g_i$ converge to their limiting value at rates greater than $\rho_{g_i} - \varepsilon$.
Then by what is proved earlier, the sum $g_1(t) + g_2(t)$ converges to $g_1^* + g_2^*$ at rate that is at least $\min\{\rho_{g_1}, \rho_{g_2}\} - \varepsilon$. But since $\varepsilon > 0$ is arbitrary, the rate of convergence is at least $\min\{\rho_{g_1}, \rho_{g_2}\}$. Similar considerations for the other two cases: scalar product and product of functions, completes the proof of the desired result.
\end{proof}

\begin{lemma} \label{lemma:reciprocal}
Let $g:\R_{\ge 0} \to \R$ be a real-valued function that converges to a nonzero real constant $g^* \in \R \setminus \{0\}$ at a rate that is at least $\rho_{g}$.  
Then $1/g(t)$ converges to $1/g^*$ at rate that is at least $\rho_g$. 
\end{lemma}
\begin{proof}
As in the proof of Lemma \ref{lem:sumproduct}, we first assume that the convergence rate of $g$ is greater than $\rho_g$.
As in the previous lemma, 
\[
 \abs{g^*} - \abs{g(t)}  \le \abs{\abs{g(t)} - \abs{g^*}} \le \abs{g(t) - g^*} \le c e^{-\rho_{g} t}, 
\]
for some positive constant $c > 0$, and so 
\[
\abs{g(t)} \ge \abs{g^*} - c e^{-\rho_{g} t}. 
\]
In particular, there is a finite time $T_0 \ge 0$ such that for all $t > T_0$, $\abs{g(t)} \ge \abs{g^*}/2$. 
It then follows that for any $t > T_0$, 
\begin{align*}
\abs{\frac{1}{g(t)} - \frac{1}{g^*}} = \abs{\frac{g(t) - g^*}{g(t) g^*}} \le \frac{2c e^{-\rho_{g} t}}{(g^*)^2 }. 
\end{align*}
But then, and as in the proof of Lemma \ref{lem:ratespeed},  we can choose a $c' > 0$ sufficiently large so that for all $t \ge 0$, 
\begin{align*}
\abs{\frac{1}{g(t)} - \frac{1}{g^*}}  \le c' e^{-\rho_{g} t}. 
\end{align*}
Finally, arguing similarly to the proof of Lemma \ref{lem:sumproduct}, we get the conclusion for the case when the convergence rate of $g$ is exactly $\rho_g$.
\end{proof}

\begin{lemma}
\label{lemma:mthroot}
Let $g:\R_{\ge 0} \to \R$ be a real-valued function that converges to a non-negative constant $g^* \in \R_{\ge 0}$ at a rate that is at least $\rho_{g}$.  
Then for any $m \in \R_{> 0}$, $g(t)^{1/m}$ converges to $(g^*)^{1/m}$ at rate that is at least
\begin{align}
\begin{cases}
\rho_{g}, &\mbox{ if } g^* > 0, \\
\rho_{g}/m, &\mbox{ if } g^* = 0. 
\end{cases}
\end{align} 
\end{lemma}
\begin{proof}
Similar to the proof of Lemma \ref{lem:sumproduct}, it suffices to assume that the convergence rate of $g$ is greater than $\rho_g$.
By Lemma \ref{lem:ratespeed} there exists a $c > 0$ such that 
\[
\abs{g(t) - g^*} < ce^{-\rho_g t}.
\]
Suppose first that $g^* > 0$. Then there exists a $T_0 \ge 0$ such that $g(t) > g^*/2$ for all $t > T_0$. 
Then for $t > T_0$, we have 
\begin{align*}
\abs{g(t) - g^*} &= \abs{g(t)^{1/m} - (g^*)^{1/m}} \abs{\sum_{i=0}^{m-1}  g(t)^{i}  (g^*)^{m-i-1}}. \\ 
&\ge \abs{g(t)^{1/m} - (g^*)^{1/m}} \sum_{i=0}^{m-1} \frac{1}{2^i} (g^*)^{i}  (g^*)^{m-i-1} \\
&= \abs{g(t)^{1/m} - (g^*)^{1/m}} \left(2^m - 1\right)\left(\frac{g^*}{2}\right)^{m-1}. 
\end{align*} 
So for $t > T_0$
\begin{align*}
\abs{g(t)^{1/m} - (g^*)^{1/m}} 
\le \left(\frac{2}{g^*}\right)^{m-1} \frac{1}{2^m-1} ce^{-\rho_g t}, 
\end{align*}
and, again by reasoning similar to the end of the proof of Lemma \ref{lem:ratespeed}, $g(t)^{1/m}$ converges to $(g^*)^{1/m} > 0$ at rate that is at least $\rho_g$. 

When $g^* = 0$, there is a $c>0$ such that 
\[
\abs{g(t)} < ce^{-\rho_g t}, 
\]
and so taking the $m$th root on both sides gives the result immediately. 
\end{proof}

\subsection{Analysis of certain non-autonomous systems}

Here we present our main mathematical results.  These will be applied to our chemical systems in section \ref{sec:elem_comp}.  

\begin{lemma}
\label{thm:system1_analysis}
Let $g_1:\R_{\ge 0} \to \R$ be a real-valued function that converges to $g_1^* \in \R$ at a rate that is at least $\rho_{g_1}$; let $g_2:\R_{\ge 0} \to \R$ be a real-valued function that converges to a positive limit $g_2^* \in \R_{> 0}$ at a rate that is at least $\rho_{g_2}$. 
We assume that the $g_1, g_2$ are smooth enough so that for any $x(0) = x_0 \ge 0$, the following non-autonomous differential equation has a unique solution $x: \R_{\ge 0} \to \R$ for all time 
\begin{align}\label{eq:linear_simple}
\dot x(t) &= g_1(t) - g_2(t)x(t). 
\end{align}
Then $x(t)$ converges to $g_1^*/g_2^*$ at rate that is at least
\begin{align}
\min\{\rho_{g_1}, \rho_{g_2}, g_2^*\}. 
\end{align}
\end{lemma}
\begin{proof}
As in the proof of Lemma \ref{lem:sumproduct}, we will assume that the convergence rates of $g_1$ and $g_2$ are greater than $\rho_{g_1}$ and $\rho_{g_2}$ respectively, and then just as in the proof of that lemma,  the result will follow immediately for the case of equality.
We first show that $x(t)$ is bounded above for all time $t \ge 0$. 
Since $g_2(t) \to g_2^* > 0$, there is a  $T_0 > 0$ such that $g_2(t) > g_2^*/2$ for all $t > T_0$. 
Moreover, note that there must be a positive constant $c_1$ such that for all time $t > 0$, we have 
\[
g_1(t) \le g_1^* + c_1 e^{-\rho_{g_1} t} \le g_1^* + c_1. 
\]
Combining the two bounds gives an upper bound for $x(t)$. Since for any $t > T_0$, 
\begin{align*}
\dot x(t)  = g_1(t) - g_2(t) x(t) \le g_1^* + c_1 - \frac12 g_2^*\cdot x(t), 
\end{align*}
we have that for any $t \ge 0$
\begin{align} \label{eq:upperboundx}
x(t) \le \max\left\{x\left([0,T_0]\right),\frac{2(g_1^*+c_1)}{g_2^*}\right\}. 
\end{align}

Now we show that $x(t)$ has the claimed rate of convergence to its limiting value $g_1^*/g_2^*$. 
First consider the differential equation obtained by replacing $g_2(t)$ with its limiting value $g_2^*$ 
     \begin{equation}\label{eq:almost_autonomous1}
         \dot z(t) = g_1(t) - g_2^* \cdot z(t), 
    \end{equation}
which has the unique solution 
\begin{equation*}
z(t) = e^{-g_2^*  t} \left( z(0)  + \int_0^t e^{g_2^*s} g_1(s) ds \right). 
\end{equation*}
Therefore,
\begin{align*}
    z(t)-\frac{g_1^*}{g_2^*} &= e^{-g_2^*  t} \left( z(0)  + \int_0^t e^{g_2^*s} g_1(s) ds \right) -\frac{g_1^*}{g_2^*}\\
    &=e^{-g_2^*  t} \left( z(0)-\frac{g_1^*}{g_2^*} + \int_0^t e^{g_2^*s} \left(g_1(s)-g_1^*\right) ds \right), 
\end{align*}
and so 
\begin{align}
\abs{z(t) - \frac{g_1^*}{g_2^*}}
&\le  e^{-g_2^*  t} \left( \left|z(0)-\frac{g_1^*}{g_2^*}\right|   +  \int_0^t  e^{g_2^* s} \abs{g_1(s)-g_1^*} ds \right) \nonumber \\
&\le  e^{-g_2^*  t} \left( \left|z(0)-\frac{g_1^*}{g_2^*}\right|   +  c \int_0^t  e^{\left(g_2^*-\rho_{g_1}\right) s}  ds \right), 
\end{align}
for some $c > 0$. 
The integral on the right hand side evaluates to 
\begin{align}
\int_0^t  e^{\left(g_2^*-\rho_{g_1}\right) s}  ds 
=
\begin{cases}
\left(e^{(g_2^*-\rho_{g_1}) t } - 1\right)/(g_2^*-\rho_{g_1}), &\mbox{ if } g_2^*\ne \rho_{g_1} \\
t,  &\mbox{ if } g_2^* = \rho_{g_1}, 
\end{cases}
\end{align}
and so 
\begin{align}
\abs{z(t) - \frac{g_1^*}{g_2^*}}
\le  
\begin{cases}
\ds e^{-g_2^*  t} \left|z(0)-\frac{g_1^*}{g_2^*}\right|   +  \frac{c}{g_2^*-\rho_{g_1}} \left(e^{-\rho_{g_1} t } - e^{-g_2^*  t}\right) , &\mbox{ if } g_2^* \ne \rho_{g_1}, \\
\ds e^{-g_2^*  t}  \left|z(0)-\frac{g_1^*}{g_2^*}\right|   +  c t e^{-g_2^*  t} , &\mbox{ if } g_2^* = \rho_{g_1}. 
\end{cases}
\end{align}

After taking logarithms, then limit superiors and applying Definition \ref{def:ratespeed}, 
in either case, we have that $z(t)$ converges to $g_1^*/g_2^*$ at least at rate $\min\{g_2^*,\rho_{g_1}\}$.

Now suppose that $x(t)$ is the solution to \eqref{eq:linear_simple} with $x(0) = x_0$ and $z(t)$ is the solution to \eqref{eq:almost_autonomous1} with $z(0) = x_0$.  
We have
\begin{align*}
    \frac{d}{dt} (z(t) - x(t)) &= -g_2^*z(t) + g_2(t) x(t)\\
    &= -g_2^*(z(t) - x(t)) + (g_2(t) - g_2^*) x(t).
\end{align*}
Since $x(t)$ is bounded above, see \eqref{eq:upperboundx}, we have that $\alpha(t) \coloneqq (g_2(t)-g_2^*)x(t)$ converges to $0$ at least at rate $\rho_{g_2}$. 
Let $h(t) \coloneqq z(t) - x(t)$, and so $h(0) = 0$ and $h(t)$ satisfies the differential equation
\begin{align} \label{eq:heqn}
\dot h(t) = \alpha(t) - g_2^* h(t). 
\end{align} 
Comparing \eqref{eq:heqn} with \eqref{eq:almost_autonomous1}, $h(t) = z(t) -x(t)$ converges to $0$ at least at rate $\min\{g_2^*, \rho_{g_2}\}$. 
Finally, since $x(t) = z(t) - h(t)$, a simple use of the triangle inequality gives us that $x(t) \to g_1^*/g_2^*$ at rate that is at least $\min\{\rho_{g_1}, \rho_{g_2}, g_2^*\}$. 
\end{proof}

\begin{lemma}
\label{thm:system2_analysis}
For $i \in \{1,2\}$, let $g_i:\R_{\ge 0} \to \R$ be a real-valued function that converges to a positive limit $g_i^* \in \R_{> 0}$ at a rate that is at least $\rho_{g_i}$. 
We assume that the $g_1, g_2$ are smooth enough so that for any $x(0) = x_0 > 0$, the following non-autonomous differential equation has a unique solution $x: \R_{\ge 0} \to \R_{\ge 0}$ for all time 
\begin{align}\label{eq:nonlinear_simple}
\dot x(t) &= x(t)(g_1(t) - g_2(t)x(t)^m), \quad \quad (m \in \Z_{> 0}). 
\end{align}
Then $x(t)$ converges to $\left(g_1^*/g_2^*\right)^{1/m}$ at rate that is at least
\begin{align}
\min\{\rho_{g_1}, \rho_{g_2}, mg_1^*\}. 
\end{align}
\end{lemma}

\begin{proof}
Similar to Lemma \ref{lem:sumproduct} and Lemma \ref{thm:system1_analysis}, it suffices to assume that the convergence rate of $g_i$ is greater than $\rho_{g_i}$ for $i \in \{1,2\}$. 
Note that the positive $x$-axis, $x>0$, is invariant, i.e. if $x(0) > 0$, then $x(t) > 0$ for all $t>0$. 
In fact, we argue that $x(t)$ has a positive lower bound as well as a positive upper bound for all time $t \ge 0$. 
The reasoning is similar to the one used in the previous lemma. 
Let $i \in \{1,2\}$. Since $g_i(t) \to g_i^* > 0$, there is a  $T_i > 0$ such that $g_i(t) > g_i^*/2$ for all $t > T_i$. Moreover, note that there must be a positive constant $c_i$ such that for all time $t > 0$, we have 
\[
g_i(t) \le g_i^* + c_i e^{-\rho_{g_i} t} \le g_i^* + c_i. 
\]
For the upper bound of $x(t)$, consider that  
\begin{align*}
\dot x(t)  = x(t) \left(g_1(t) - g_2(t) x(t)^m \right) \le x(t) \left( g_1^* + c_1 - \frac12 g_2^*\cdot x(t)^m \right), 
\end{align*}
for any $t > T_2$. Then we have that for any $t > 0$
\begin{align} \label{eq:upperboundx2}
x(t) \le \max\left\{x\left([0,T_2]\right), \left( \frac{2(g_1^*+c_1)}{g_2^*} \right)^{1/m}\right\}. 
\end{align}
For the positive lower bound of $x(t)$, consider that  
\begin{align*}
\dot x(t)  = x(t) \left(g_1(t) - g_2(t) x(t)^m \right) \ge x(t) \left( g_1^*/2 -  (g_2^* + c_2) \cdot x(t)^m \right), 
\end{align*}
for any $t > T_1$. Then we have that for any $t > 0$
\begin{align} \label{eq:lowerboundx2}
x(t) \ge \min\left\{x\left([0,T_1]\right), \left( \frac{g_1^*}{2(g_2^*+c_2)} \right)^{1/m}\right\}. 
\end{align}

To show that we have the desired rate of convergence, we let  
 \begin{align*}
    y(t) \coloneqq \frac{1}{x(t)^m}.
 \end{align*}
Differentiating  yields 
\begin{align*}
    \dot y(t) = - \frac{m}{x(t)^{m+1}}\cdot  \dot x(t)= mg_2(t) - m g_1(t) y(t). 
\end{align*}
The functions $mg_2$ and $mg_1$ satisfy the hypotheses of Lemma \ref{thm:system1_analysis}, and so we have that $y(t)$ converges to $g_2^*/g_1^*$ at rate that is at least $\min\{\rho_{g_1},\rho_{g_2},mg_1^*\}$. 
Then applying Lemma \ref{lemma:reciprocal} gives us that the reciprocal of $y(t)$, i.e. $x(t)^m$ converges to $g_1^*/g_2^*$ at rate that is at least $\min\{\rho_{g_1},\rho_{g_2},mg_1^*\}$. 
Finally applying Lemma \ref{lemma:mthroot} gives us that $x(t)$ converges to $\left(g_1^*/g_2^*\right)^{1/m}$ at rate that is at least $\min\{\rho_{g_1},\rho_{g_2},mg_1^*\}$. 
\end{proof}

\begin{lemma}
\label{thm:test2}
Let $g_1:\R_{\ge 0} \to \R$ be a real-valued function that converges to a negative limit $g_1^* \in \R_{< 0}$ at a rate that is at least $\rho_{g_1}$; let $g_2:\R_{\ge 0} \to \R$ be a real-valued function that converges to a positive limit $g_2^* \in \R_{> 0}$ at a rate that is at least $\rho_{g_2}$. 
We assume that  $g_1$ and $g_2$ are smooth enough so that for any $x(0) = x_0 > 0$, the following non-autonomous differential equation has a unique solution $x: \R_{\ge 0} \to \R_{\ge 0}$ for all time 
\begin{align}\label{eq:nonlinear_simple_np}
\dot x(t) &= x(t)(g_1(t) - g_2(t)x(t)^m), \quad \quad (m \in \Z_{> 0}). 
\end{align}
Then $x(t)$ converges to $0$ at rate that is at least
\begin{align}
\min\{\rho_{g_1}, - g_1^*\}. 
\end{align}
\end{lemma}
\begin{proof}
Similar to the previous results,  it suffices to assume that the convergence rate of $g_i$ is greater than $\rho_{g_i}$ for $i \in \{1,2\}$.
Since $g_2(t) \to g_2^* > 0$, there is a  $T_0 > 0$ such that $g_2(t) > g_2^*/2$ for all $t \ge T_0$. 
So for any $t \ge T_0$, we have that 
\begin{align*}
\dot x(t) \le g_1(t)x(t), 
\end{align*}
which implies that $0 \le x(t) \le z(t)$ where $z(t)$ satisfies $z(T_0) = x(T_0)$ and for any $t \ge T_0$, 
\begin{align*}
\dot z(t) = g_1(t)z(t). 
\end{align*}
This last equation satisfies the hypotheses of Lemma \ref{thm:system1_analysis} with $g_1(t)$ of the lemma taken to be identically $0$ and $g_2(t)$ of the lemma replaced with $-g_1(t)$. The conclusion is that $z(t)$ and therefore $x(t)$ converges to $0$ at rate that is at least $\min\{\rho_{g_1}, - g_1^*\}$. 
\end{proof}

\begin{lemma}
\label{thm:test}
Let $g:\R_{\ge 0} \to \R$ be a real-valued function that converges to $0$ at a rate that is at least $\rho_{g}$. 
We assume that $g$ is smooth enough so that for any $x(0) = x_0 > 0$, the following non-autonomous differential equation has a unique solution $x: \R_{\ge 0} \to \R_{\ge 0}$ for all time 
    \begin{align*}
        \dot x(t) = x(t) (1 - g(t) x(t)^m), \quad \quad (m \in \Z_{> 0}). 
    \end{align*}
Then for any $\varepsilon > 0$, there is a $c_\varepsilon >0$ such that
    \[
    x(t) \ge c_\varepsilon e^{(\min\{\rho_g/m,1\} - \varepsilon)t}, \text{ for all } t \ge 0.
    \]
\end{lemma}

\begin{proof}
Note that the positive $x$-axis is invariant. We may therefore define
 \begin{align*}
    y(t) \coloneqq \frac{1}{x(t)^m},   
 \end{align*}
which satisfies the differential equation 
\begin{align*}
    \dot y(t) = m g(t) - m y(t).
\end{align*}
Since $y(t)$ satisfies the hypotheses of Lemma \ref{thm:system1_analysis} with $g_1(t) \coloneqq m g(t)$ and $g_2(t) \equiv m$, we have that $y(t)$ converges to $0$ at rate that is at least $\min\{\rho_g,m\}$. Hence, for every $\varepsilon > 0$ there is a positive constant $c_\varepsilon$ such that 
\[
y(t) \le c_\varepsilon e^{- (\min\{\rho_g,m\} - m\cdot \varepsilon) t}. 
\]
In terms of the original variable $x(t)$ we have that 
\[
\frac{1}{x(t)^m} \le c_\varepsilon e^{-( \min\{\rho_g,m\} - m\cdot \varepsilon) t}. 
\]
Taking the reciprocal and then the $m$th root gives the desired inequality.  
\end{proof}

\section{Elementary operations} \label{sec:elem_comp}

We construct reaction network  modules that can carry out the elementary arithmetic operations of: identification, inversion, $m$th root (for $m \ge 2$), addition, multiplication, absolute difference, rectified subtraction, and partial real inversion (for reals represented via the dual rail representation).  These operations form the core of basic arithmetic. 
 Compositions of them will be taken up in the next section. 
Each of the elementary operations that we consider has one or two non-negative real inputs, hence referred to as a unary or binary operation, and one non-negative real output. 
In this section, we will show that the performance of each elementary operation retains from the input function the desired features of global exponential convergence, and convergence speed that is bounded from below independent of the inputs and output of computation. 
We will apply results from the previous section to prove these properties. 

Throughout this section all rate constants are assumed to be equal to 1.  Moreover, we assume throughout that the input functions (typically denoted as $a(t)$ or $b(t)$) are assumed to be smooth enough so that the associated differential equations have unique solutions for all relevant initial conditions and all time.

\subsection{Identification}

\begin{lemma}[Identification] \label{thm:identification}
Consider the reaction network and mass-action system
\begin{align}
\{A \xrightarrow{} A + X, \quad X \xrightarrow{} 0\}, \quad \quad \dot x(t) = a(t) - x(t),  
\end{align}
where $a(t)$ is a non-negative-valued function of time that converges to a non-negative constant $a^*$ at a rate that is at least $\rho_a$. 
Then the concentration of species $X$, i.e. the variable $x$, computes identification $a^* \mapsto x^* = a^*$ at speed that is at least $\min\{\rho_{a},1\}$.
\end{lemma}
\begin{proof}
The result follows immediately by applying Lemma \ref{thm:system1_analysis} with  $g_1(t) \coloneqq a(t)$ and $g_2(t) \equiv 1$. 
\end{proof}

 Before proceeding further, we present an example of a composition of two identifications, which is arguably the simplest possible composition.  We do so to point out how to utilize our results in that context.  This topic will be extensively revisited in section \ref{sec:combining}.
\begin{example} 
\label{ex:twoids}
    Consider a sequence of two identification operations. We will use $A$ to denote the input species, whose concentration is unchanged over time. Let $Y$ denote a species whose long-term concentration is identified with that of $A$ and $X$ the species whose long-term value is identified with the long-term value of $Y$. The reaction network and system of ODEs for this computation are: 
    \begin{align*}
        A \to A + Y, \quad &Y \to 0, \quad Y \to Y + X, \quad X \to 0. \\
        \dot y(t) &= a - y(t) \\
        \dot x(t) &= y(t) - x(t). 
    \end{align*}
    We can solve this system explicitly and thereby obtain exact rates of convergence. The solution to the initial value problem with $x(0) = x_0$ and $y(0) = y_0$ is 
        \begin{align*}
        y(t) &= a\left(1 - e^{-t}\right) + y_0 e^{-t} \\ 
        x(t) &= a\left(1 - (1+t)e^{-t}\right) + (x_0 + ty_0) e^{-t}
    \end{align*}
Thus both $y(t)$ and $x(t)$ converge to $a$ at speed equal to $1$. 
This is, of course, consistent with the conclusion from Lemma \ref{thm:identification}. Since the input species $A$ is not changing, its rate of convergence is $\rho_a = + \infty$. By Lemma \ref{thm:identification}, $y(t)$ converges at speed $\rho_y$ that is at least $\min \{\rho_a, 1\} = \min \{\infty, 1\} = 1$. 
Again, by application of Lemma \ref{thm:identification}, $x(t)$ converges at speed that is at least $\min \{\rho_y, 1\} = \min \{1, 1\} = 1$. 
\end{example}

\subsection{Inversion}

\begin{lemma}[Inversion] \label{thm:inversion}
Consider the reaction network and mass-action system
\begin{align}
\{X \to 2X, \quad A + 2X \to A + X\}, \quad \quad \dot x(t) = x(t)(1 -a(t)x(t)), 
\end{align}
where $a(t)$ is a non-negative-valued function of time that converges to a non-negative constant $a^*$ at a rate that is at least $\rho_a$. 
Then the concentration of species $X$, i.e. the variable $x$, computes inversion $a^* \mapsto x^* = 1/a^*$ at speed that is at least $\min\{\rho_{a},1\}$.
\end{lemma}
\begin{proof}
The result follows immediately Lemma \ref{thm:system2_analysis} with  $m=1$, $g_1(t) \equiv 1$, and $g_2(t) \coloneqq a(t)$. 
\end{proof}

\subsection{$m$th root}

\begin{lemma}[$m$th root] \label{thm:mthroot}
Consider the reaction network and mass-action system
\begin{equation} \label{eq:mthroot}
\begin{aligned}
\{Y \to 2Y, \quad A + (m+1)Y &\to A + mY, \quad X \to 2X, \quad Y + 2X \to Y + X\} \\
\dot y(t) &= y(t) \left(1 - a(t)y(t)^m\right), \\
\dot x(t) &= x(t) \left(1 - y(t)x(t)\right),
\end{aligned}
\end{equation}
where $m \in \Z_{\ge 2}$ and $a(t)$ is a non-negative-valued function of time that converges to a non-negative constant $a^*$ at a rate that is at least $\rho_a$. 
Then the concentration of species $X$, i.e. the variable $x$, computes the $m$th root $a^* \mapsto x^* = \sqrt[m]{a^*}$ at speed that is at least
\begin{align} \label{eq:speedmthroot} 
\begin{cases}
\min\{\rho_a, 1\}, &\mbox{ if } a^* \ne 0,  \\
\min\left\{\frac{\rho_a}{m}, 1 \right\}, &\mbox{ if } a^* = 0. 
\end{cases}
\end{align}  
\end{lemma}
\begin{proof}
We first consider the case $a^* > 0$.
The first differential equation satisfies the hypotheses of  Lemma \ref{thm:system2_analysis} with $m \ge 2$, $g_1(t) \equiv 1$ and $g_2(t) \coloneqq a(t)$, and therefore $y(t)$ converges to $(a^*)^{-1/m}$  at rate that is at least $\min\{\rho_{a}, m\}$. 
Application of the same lemma again to the second differential equation with $m =1$, $g_1(t) \equiv 1$ and $g_2(t) \coloneqq y(t)$ gives us that $x(t)$ converges to $\sqrt[m]{a^*}$  at rate that is at least $\min\{\rho_{a}, 1\}$. 

Now we consider the case when $a^*=0$. We apply Lemma \ref{thm:test} to the first differential equation with $m \ge 2$ and $g(t) \coloneqq a(t)$.  
Fix $\varepsilon > 0$. Then by Lemma \ref{thm:test}, there exists a constant $c_\varepsilon > 0$ such that 
\[
y(t) \ge c_\varepsilon e^{\alpha t}, 
\]
where $\alpha \coloneqq \min\{\rho_a/m, 1\} - \varepsilon$.
Therefore,
\begin{align} \label{eq:temp1}
    \dot x(t) \le x(t) \left(1 - c_\varepsilon e^{\alpha t}x(t)\right).
\end{align}
Define a change of variables 
\[
z(t) = e^{\alpha t} x(t).
\]
Clearly $z(0) = x(0)$ and $\dot z(t) = \alpha z(t) + e^{\alpha t} \dot x(t)$. 
Multiplying both sides of \eqref{eq:temp1} with $e^{\alpha t}$, we get 
\begin{align*}
e^{\alpha t} \dot x(t) \le e^{\alpha t}x(t) \left(1 - c_\varepsilon e^{\alpha t}x(t)\right), 
\end{align*}
which we may rewrite as 
\begin{align*}
    \dot z(t) \le z(t) \left(1 + \alpha - c_\varepsilon z(t)\right).
\end{align*}
The solution $z(t)$ is bounded above by the solution of the autonomous system 
\begin{align*}
    \dot w(t) = w(t) \left(1 + \alpha - c_\varepsilon w(t)\right), \quad w(0) = z(0). 
\end{align*}
The explicit solution gives an upper bound for $x(t)$ 
\[
x(t) = e^{-\alpha t} z(t) \le e^{-\alpha t} w(t) = \frac{1 + \alpha}{c_\varepsilon(e^{\alpha t} - e^{-t}) + \frac{1+\alpha}{x_0} e^{-t}}. 
\]
The right hand side converges to $0$ at a rate that is at least $\alpha = \min\{\rho_a/m,1\} - \varepsilon$. 
But since the $\varepsilon > 0$ is arbitrary, the convergence is at rate $\min\{\rho_a/m,1\}$. 
\end{proof}

\subsection{Addition}

\begin{lemma}[Addition] \label{thm:addition}
Consider the reaction network and mass-action system
\begin{align}
\{A \to A + X, \quad B \to B + X, \quad X \to 0\}, \quad \quad \dot x(t) = a(t) + b(t) - x(t), 
\end{align}
where $a(t)$ and $b(t)$ are non-negative functions of time that converge to non-negative constants $a^*$ and $b^*$ at rates that are at least $\rho_a$ and $\rho_b$, respectively. 
Then the concentration of species $X$, i.e. the variable $x$, computes addition $(a^*,b^*) \mapsto x^* = a^* + b^*$ at speed that is at least $\min\{\rho_{a}, \rho_b, 1\}$.
\end{lemma}
\begin{proof}
By Lemma \ref{lem:sumproduct}, the sum $a(t) + b(t)$ converges to $a^*+b^*$ at speed that is at least $\min\{\rho_a, \rho_b\}$. 
The result then follows from Lemma \ref{thm:system1_analysis} with $g_1(t) \coloneqq a(t) + b(t) $ and $g_2(t) \equiv 1$. 
\end{proof}

\subsection{Multiplication}

\begin{lemma}[Multiplication] \label{thm:multiplication}
Consider the reaction network and mass-action system
\begin{align}
\{A + B \to A + B + X, \quad X \to 0\}, \quad \quad \dot x(t) = a(t)b(t) - x(t), 
\end{align}
where $a(t)$ and $b(t)$ are non-negative functions of time that converge to non-negative constants $a^*$ and $b^*$ at rates that are at least $\rho_a$ and $\rho_b$, respectively. 
Then the concentration of species $X$, i.e. the variable $x$, computes multiplication $(a^*,b^*) \mapsto x^* = a^* \cdot b^*$ at speed that is at least
\begin{align}
\begin{cases}
\min\{\rho_{a}, \rho_{b}, 1\}, &\mbox{ if } a^* \ne 0,  b^* \ne 0, \\
\min\{\rho_{a}, 1\}, &\mbox{ if } a^* = 0,  b^* \ne 0, \\
\min\{\rho_{b}, 1\}, &\mbox{ if } a^*  \ne 0,  b^* = 0, \\
\min\{\rho_{a} + \rho_{b}, 1\}, &\mbox{ if } a^* = 0,  b^* = 0. 
\end{cases}
\end{align}
\end{lemma}
\begin{proof}
Applying Lemma \ref{lem:sumproduct} to the product  $a(t) b(t)$ gives the lower bound on its rate of convergence. 
The result then follows by applying Lemma \ref{thm:system1_analysis} with $g_1(t) \coloneqq a(t)  b(t) $ and $g_2(t) \equiv 1$. 
\end{proof}

\subsection{Absolute difference}

\begin{lemma}[Absolute difference] \label{thm:absdiff}
Consider the reaction network and mass-action system
\begin{equation} \label{eq:absdiff}
\begin{aligned}
\{Y \to 2Y, \quad 2A + 3Y \to 2A + 2Y, &\quad 2B + 3Y \to 2B + 2Y, \quad A + B + 3Y \to A + B + 5Y, \\ 
&X \to 2X, \quad Y + 2X \to Y + X\}, \\
\dot y(t) &= y(t) \left(1 - (a(t)-b(t))^2 y(t)^2\right), \\
\dot x(t) &= x(t) \left(1 - y(t)x(t)\right),
\end{aligned}
\end{equation}
where $a(t)$ and $b(t)$ are non-negative functions of time that converge to non-negative constants $a^*$ and $b^*$ at rates that are at least $\rho_a$ and $\rho_b$, respectively.  
Then the concentration of species $X$, i.e. the variable $x$, computes absolute difference $(a^*,b^*) \mapsto x^* = \abs{a^* - b^*}$ at speed that is at least
\begin{align} \label{eq:speedabsdiff}
\min\{\rho_a, \rho_b, 1\}. 
\end{align}  
\end{lemma}
\begin{proof}
By Lemma \ref{lem:sumproduct}, the difference $a(t) - b(t)$ converges to $a^* - b^*$ at speed that is at least $\min\{\rho_a, \rho_b\}$. 
An application of the product part of the same Lemma \ref{lem:sumproduct} gives that $(a(t) - b(t))^2$ converges to $(a^* - b^*)^2$ at rate that is at least
\[
\begin{cases}
\min\{\rho_a, \rho_b\}, &\mbox{ if } a^* \ne b^*\\
2\min\{\rho_a, \rho_b\}, &\mbox{ if } a^* = b^*. 
\end{cases}
\]
Finally, we apply Lemma \ref{thm:mthroot}  to the system \eqref{eq:absdiff} with $m=2$ and the function $a(t)$ replaced with $(a(t) -b(t))^2$. 
\end{proof}

\subsection{Rectified subtraction}

\begin{lemma}[Rectified subtraction]
Consider the reaction network and mass-action system
\begin{equation} \label{eq:rectifiedsubtraction}
\begin{aligned}
\{Y \to 2Y, \quad 2A + 3Y \to 2A + 2Y, &\quad 2B + 3Y \to 2B + 2Y, \quad A + B + 3Y \to A + B + 5Y, \\ 
A + Y + X  \to A + Y + 2X, \quad &B + Y  + X \to B + Y, \quad Y + 2X \to Y + X\}, \\
\dot y(t) &= y(t) \left(1- (a(t)-b(t))^2y(t)^2\right),  \\
\dot x(t) &= y(t) x(t) \left( a(t) - b(t) - x(t) \right), 
\end{aligned}
\end{equation}
where $a(t)$ and $b(t)$ are non-negative functions of time that converge to non-negative constants $a^*$ and $b^*$ at rates that are at least $\rho_a$ and $\rho_b$, respectively. 
Then the concentration of species $X$, i.e. the variable $x$, computes rectified subtraction 
\[
(a^*,b^*) \mapsto 
\begin{cases} 
a^*-b^* &\mbox{ if } a^* > b^* \\ 
0 &\mbox{ if } a^* \le b^* 
\end{cases} 
\]
at speed that is at least
\begin{align} \label{eq:speedrectsub}
\min\{\rho_a, \rho_b, 1\}. 
\end{align} 
\end{lemma}
\begin{proof}
By Lemma \ref{lem:sumproduct}, the difference $a(t) - b(t)$ converges to $a^* - b^*$ at rate that is at least $\rho \coloneqq \min\{\rho_a, \rho_b\}$. 

First suppose that  $a^* \ne b^*$. 
Similar to the proof of Lemma \ref{thm:absdiff}, $(a(t) - b(t))^2$ converges to $(a^* - b^*)^2$ at rate that is at least $\min\{\rho_a, \rho_b\}$. 
The differential equation for $y(t)$ satisfies the hypotheses of  Lemma \ref{thm:system2_analysis} with $m = 2$, $g_1(t) \equiv 1$ and $g_2(t) \coloneqq (a(t) - b(t))^2$, and therefore $y(t)$ converges to $\abs{a^*-b^*}^{-1} > 0$ at rate that is at least $\min\{\rho_{a}, \rho_b, 2\}$. 
Now rewrite the differential equation for $x(t)$ as 
\begin{align} \label{eq:int2nd}
\dot x(t) = x(t) \left(y(t) (a(t) - b(t)) - y(t) x(t) \right), 
\end{align}
and note that by Lemma \ref{lem:sumproduct}, $y(t) (a(t) -b(t))$ converges to either $+1$ or $-1$, depending on whether $a^* > b^*$ or $a^* < b^*$, at rate that is at least $\min\{\rho_{a}, \rho_b, 2\}$. 
When $a^* > b^*$, we apply Lemma \ref{thm:system2_analysis} with $m=1$, $g_1(t) \coloneqq  y(t) (a(t) - b(t))$, and $g_2(t) = y(t)$ to \eqref{eq:int2nd}, and we conclude that $x(t)$ converges to $a^* - b^*$ at rate that is at least $\min\{\rho_{a}, \rho_b, 1\}$. 
When $a^* < b^*$, we apply Lemma \ref{thm:test2} with $m=1$, $g_1(t) \coloneqq  y(t) (a(t) - b(t))$, and $g_2(t) = y(t)$ to \eqref{eq:int2nd}, and we conclude that $x(t)$ converges to $0$ at rate that is at least $\min\{\rho_{a}, \rho_b, 1\}$. In either case, we get the desired result when $a^* \ne b^*$. 

Now consider the case that $a^* = b^*$. 
We have that $a(t) - b(t)$ converges to $0$ at rate that is at least $\rho = \min\{\rho_a,\rho_b\}$. We fix an  $\varepsilon > 0$. Then, there is a constant $c_\varepsilon > 0$ such that for all $t \ge 0$, 
\begin{align} \label{eq:rateaminusb}
\abs{a(t) - b(t)} < c_\varepsilon e^{-(\rho - \varepsilon) t}. 
\end{align}
Similar to the proof of Lemma \ref{thm:absdiff}, $(a(t) - b(t))^2$ converges to $0$ at rate that is at least $2\min\{\rho_a, \rho_b\}$. 
The differential equation for $y(t)$ satisfies the hypotheses of  Lemma \ref{thm:test} with $m = 2$, and $g(t) \coloneqq (a(t) - b(t))^2$, and therefore there is a $c_\varepsilon'>0$ such that for all $t \ge 0$, we have 
    \begin{align} \label{eq:rateygrowth}
    y(t) \ge c_\varepsilon' e^{(\min\{\rho_a, \rho_b, 1\}-\varepsilon) t} = c_\varepsilon' e^{\nu t},
    \end{align}
    where $\nu \coloneqq \min\{\rho_a, \rho_b, 1\} - \varepsilon$. 

We will now show that $x(t)$ converges to $0$ at rate that is at least $\nu$.  As in our previous results, this will be enough to prove the result as $\varepsilon>0$ is arbitrary.
Note that the positive $x$-axis is invariant, i.e. for any $x(0) > 0$, $x(t) > 0$ for all time $t \ge 0$. 
Define 
\[
u(t) \coloneqq \frac{1}{K} x(t) e^{\nu t}, \quad \mbox{ where } K \coloneqq c_\varepsilon + (\nu+1)/c_\varepsilon', 
\]
and so $u(t) > 0$ for all $t \ge 0$ whenever $x(0) > 0$. 
It suffices to show that if $u(t) \ge 1$, then $\dot u(t) \le -1$. 
Indeed assuming this is true, we first show that $x(t)$ converges at the desired rate. Consider two cases: (i) $u(0) < 1$, since the region $u < 1$ is invariant we have that 
\[
x(t) = K u(t) e^{-\nu t} \le K e^{-\nu t},
\]
and (ii) $u(0) \ge 1$, in which case we have that $u(t) \le u(0)$ for all $t \ge 0$, and so 
\[
x(t) = K u(t) e^{-\nu t} \le K u(0) e^{-\nu t} = x(0) e^{-\nu t}. 
\]
Combining the two cases, we have that for all $t \ge 0$,
\[
x(t) \le \max\{x(0), K\} e^{-\nu t},
\]
i.e. $x$ converges to $0$ at rate that is at least $\nu$.
Since $\varepsilon > 0$ is arbitrary, $x$ converges to $0$ at rate that is at least $\min\{\rho_a, \rho_b, 1\}$.

Now we prove the claim that if $u(t) \ge 1$ then $\dot u(t) \le -1$. 
Note that we have 
\begin{align*}
\dot u &= \frac{e^{\nu t}}{K} \left( \dot x + \nu x \right) = \nu u + u \frac{\dot x}{x} \\
&= \nu u + uy(a - b - K u e^{-\nu t}) \le \nu u + uy \left(c_\varepsilon e^{(-\rho + \varepsilon) t} - K u e^{-\nu t}\right) \le \nu u + uy(c_\varepsilon - K u) e^{-\nu t}, 
\end{align*}
where we used \eqref{eq:rateaminusb}, $\nu = \min\{\rho,1\} - \varepsilon \le \rho - \varepsilon$, and also that $u(t)$ and $y(t)$ are positive for all $t \ge 0$. Since we are assuming that $u \ge 1$, we then have 
\begin{align*}
\dot u &\le \nu u + uy(c_\varepsilon - K) e^{-\nu t} = \nu u - uy \left(\frac{\nu + 1}{c_\varepsilon'} \right) e^{-\nu t} \\ 
&\le \nu u - u \left(\nu + 1\right) = -u \le -1, 
\end{align*}
where we used \eqref{eq:rateygrowth} and the assumption $u \ge 1$. This concludes the proof. 
\end{proof}

\subsection{Partial real inversion}

The input to the elementary operation of \textit{partial real inversion} is a nonzero real number and can be positive or negative.  Thus, the dual rail representation is used for the input, see section \ref{sec:negativity}.

\begin{lemma}[Partial real inversion] \label{thm:partialrealinversion}
Consider the reaction network and mass-action system
\begin{equation} \label{eq:partialrealinversion}
\begin{aligned}
\{Y \to 2Y, &\quad A_p + 2Y \to A_p + Y, \quad A_n + 2Y \to A_n + Y,  \\ 
A_p + Y + X   \to A_p + Y + 2X, \quad &A_n + Y + X   \to A_n + Y, \quad 2A_p + Y + 2X   \to 2A_p + Y + X\}, \\
\dot y(t) &= y(t) \left(1- (a_p(t) + a_n(t))y(t) \right) \\
\dot x(t) &= y(t) x(t) \left( a_p(t)(1- a_p(t)x(t)) - a_n(t)\right) 
\end{aligned}
\end{equation}
where $a_p(t)$ and $a_n(t)$ are non-negative functions of time that converge to non-negative constants $a_p^*$ and $a_n^*$ at rates that are at least $\rho_{a_p}$ and $\rho_{a_n}$, respectively. 
We assume that one and only one of $a_p^*$ or $a_n^*$ is positive, while the other  is zero. 
Then the concentration of species $X$, i.e. the variable $x$, computes partial real inversion 
\[
(a_p^*,a_n^*) \mapsto 
\begin{cases} 
1/a_p^*,  &\mbox{ if } a_p^* > 0,  \\ 
0,  &\mbox{ if } a_p^* = 0, 
\end{cases} 
\]
at speed that is at least
\begin{align} \label{eq:speedpartialrealinversion}
\min\{\rho_{a_p}, \rho_{a_n}, 1\}. 
\end{align} 
\end{lemma}
\begin{proof}
By Lemma \ref{lem:sumproduct}, the sum $a_p(t) + a_n(t)$ converges to $a_p^* + a_n^* > 0$ at rate that is at least $\min\{\rho_{a_p}, \rho_{a_n}\}$. 
The first differential equation in \eqref{eq:partialrealinversion} satisfies the hypotheses of  Lemma \ref{thm:system2_analysis} with $m = 1$, $g_1(t) \equiv 1$ and $g_2(t) \coloneqq a_p(t) + a_n(t)$, and therefore $y(t)$ converges to $(a_p^* + a_n^*)^{-1}$  at rate that is at least $\min\{\rho_{a_p}, \rho_{a_n}, 1\}$. 

The differential equation for $x(t)$ can be written as 
\[
\dot x(t) = x(t) \left( y(t) (a_p(t) - a_n(t)) - y(t) a_p(t)^2x(t) \right)
\]
By Lemma \ref{lem:sumproduct}, the difference $a_p(t) -a_n(t)$ converges to its limit $a_p^* - a_n^*$ at rate that is at least $\min\{\rho_{a_p}, \rho_{a_n}\}$.  
By the product part of the same lemma, the product $y(t) (a_p(t) -a_n(t))$ converges to $+1$ or $-1$, depending on whether $a_p^* >0$ or $a_p^* = 0$, at rate that is at least $\min\{\rho_{a_p}, \rho_{a_n}, 1\}$. 

First consider the case $a_p^* > 0$ and $a_n^* =0$.  
By Lemma \ref{lem:sumproduct}, the product $y(t) a_p(t)^2$ converges to $a_p^*$ at rate that is at least $\min\{\rho_{a_p}, \rho_{a_n}, 1\}$. 
The differential equation for $x(t)$ satisfies the hypotheses of  Lemma \ref{thm:system2_analysis} with $m = 1$, $g_1(t) \coloneqq y(t) (a_p(t) -a_n(t))$ and $g_2(t) \coloneqq y(t) a_p(t)^2$, and therefore $x(t)$ converges to $1/a_p^*$ at rate that is at least $\min\{\rho_{a_p}, \rho_{a_n}, 1\}$.

Now consider the case $a_p^* = 0$ and $a_n^* > 0$.  
Suppose that $z(t)$ satisfies the differential equation $\dot z(t) = z(t)  y(t) (a_p(t) - a_n(t))$, then by applying Lemma \ref{thm:system1_analysis} with $g_1(t) \equiv 0$, and $g_2(t) \coloneqq - y(t) (a_p(t) - a_n(t))$, we get that $z(t)$ converges to $0$ at rate that is at least $\min\{\rho_{a_p}, \rho_{a_n}, 1\}$. 
Since we have 
\begin{align*}
\dot x(t) &= x(t) \left( y(t) (a_p(t) - a_n(t)) - y(t) a_p(t)^2x(t) \right) \\
&\le x(t)  y(t) (a_p(t) - a_n(t)), 
\end{align*}
for any $x(0) = z(0)$, we must have $x(t) \le z(t)$ for all $t \ge 0$, and so $x(t)$ converges to $0$ at rate that is at least $\min\{\rho_{a_p}, \rho_{a_n}, 1\}$. 
\end{proof}

\section{Composite operations}
\label{sec:combining}

\subsection{Main theory}
\label{sec:main_theory}

The elementary operations in section \ref{sec:elem_comp} can be combined to produce composite operations. 
Some operations, such as division, maximum, subtraction, which may be considered ``elementary'' from the point of view of arithmetic are in fact composite in our construction as they are produced by combining the operations in the previous section.

We start with a theorem about the speed of a general arithmetic operation, one that involves finitely many compositions of elementary operations. 
The result presented below is the most consequential theorem of this paper.

\begin{theorem} [Composite computations] \label{thm:composite}
Consider a computation that is composed from a finite number of elementary computations from section \ref{sec:elem_comp}. 
Then the speed of this composite computation is at least that of the slowest elementary computation. 
\end{theorem}
\begin{proof}
The proof only requires combining the results in Lemmas \ref{thm:identification}-\ref{thm:partialrealinversion}. 
\end{proof}

Example \ref{ex:twoids} is helpful here. 
When doing a sequence of two identifications, the intermediate variable converges like $e^{-t}$, while the output of the second identification converges like $te^{-t}$. 
For both variables, intermediate and output, the speed of convergence is $1$. 

For all elementary operations, we used the assumption that the mass action reaction rate constants are all equal. In fact, we set the rate constants to be equal to $1$. If instead, all rate constants are taken to be $\sigma > 0$, a time change argument can be used to make them all equal to $1$ in the modified system (see Remark \ref{remark:scaling} below). We now state two corollaries which assume that the rate constants are all equal to $1$. 
\begin{corollary}
\label{cor:speed_special}
Consider a computation that is composed from a finite number of elementary computations  from section \ref{sec:elem_comp}.  
If none of the elementary computations are a root of zero, then the speed of the composite computation is at least 1.
\end{corollary}

\begin{corollary}
\label{cor:speed_general}
Consider a computation that is composed from a finite number of elementary computations  from section \ref{sec:elem_comp}.
Suppose that of all the elementary steps in  the composite computation, there are $p$ elementary steps which require taking the $m_j$-th root of zero, $1 \le j \le p$. 
Then the speed of the composite computation is at least 
\[
\prod_{j=1}^p \frac{1}{m_j}. 
\]
\end{corollary}

\begin{remark}
    \label{remark:scaling}
    Throughout this work, we have always chosen the rate constants of the reactions to be equal to 1.  However, the choice of the value 1 was arbitrary.  For example, suppose that the differential equation for our system, with a rate constant of 1, is
    \[
    \dot x(t) = f(x(t)).
    \]
    Now suppose that we instead choose rate constants of $\sigma>0$, and denote the system by
    \[
        \dot z(t) = \sigma f(z(t)).
    \]
    We then see that $x$ and $z$ are related by a simple time-change: $z(t) = x(\sigma t)$.  You can see this by simply utilizing the chain rule:
        \begin{align*}
    \frac{d}{dt} x(\sigma t) = f(x(\sigma t)) \cdot \sigma.
    \end{align*}
    Thus, we immediately conclude from Corollary \ref{cor:speed_special} that if none of the elementary computations are the root of zero, the speed of the composite computation via $z$ is $\sigma$, whereas from Corollary \ref{cor:speed_general}, we conclude the speed via $z$ is $\prod_{j=1}^p \sigma/m_j$, when the conditions of that corollary are met.  We discuss the possibility of having non-equal rate constants in section \ref{sec:discussion}.\hfill $\triangle$
\end{remark}

We now show how to compose operations from 
section \ref{sec:elem_comp}
to carry out the operations \textit{division}, \textit{maximum}, and \textit{subtraction}.

\subsubsection*{Division}

Division, $a^*/b^*$, for $a^*,b^*>0$, is a composite operation of our elementary operations.  First we compute the reciprocal $b^* \mapsto 1/b^*$, and then the product $(a^*,1/b^*) \mapsto a^*/b^*$, see Figure \ref{fig:division}. 
The reaction network for computing division can be built simply by taking the union of reaction network \textit{modules}, each of which computes an elementary operation, as follows
\begin{equation} \label{eq:division}
\begin{aligned}
\{Z \to 2Z, \quad B + 2Z \to B + Z\}, \quad \quad &\dot z(t) = z(t)(1 -b(t)z(t)) \quad \mbox{(inversion)}, \\
\{A + Z \to A + Z + X, \quad X \to 0\}, \quad \quad &\dot x(t) = a(t)z(t) - x(t) \quad \mbox{(multiplication)}. 
\end{aligned}
\end{equation}
Suppose $a(t)$ and $b(t)$ are real-valued functions that converge to $a^* \ge 0$ and $b^* > 0$, respectively. 
Then concentration of species $X$, i.e. the variable $x(t)$ computes $a^*/b^*$. 
The minimum speed of convergence can be determined as follows. 
By combining results from Lemmas  \ref{thm:inversion}, and \ref{thm:multiplication}, we have that 
\begin{align*}
\rho_x &\ge  \min\{\rho_a, \rho_z, 1\} \\
&\ge  \min\{\rho_a, \min\{\rho_b, 1\}, 1\} = \min\{\rho_a, \rho_b, 1\}. 
\end{align*}
In the previous short calculation, we have shown that specifically for a composition involving reciprocation and multiplication, the minimum speed of the composite computation is at least the minimum speed of the slowest elementary computation. 
This is, of course, an instance of the general result, Theorem \ref{thm:composite}.

\subsubsection*{Maximum}
The maximum $\max\{a,b\}$ of two positive reals $a$ and $b$ can be found using the following formula
\[
\max\{a,b\} = \frac12 \left(a + b + \abs{a-b} \right), 
\]
which is a composition of the absolute difference, two binary additions and a multiplication by $1/2$. 
By Theorem \ref{thm:composite}, the speed of computation is bounded below.

\subsubsection*{Subtraction}
The output $a-b$ of subtraction for two non-negative reals $a,b \ge 0$ is an arbitrary real number, and therefore requires the dual rail representation, i.e. two chemical species, the concentration of one is the positive part of $a-b$ and the concentration of the other is the negative part. 
We may write $a-b$ in dual rail representation as 
\[
a-b = (a \mathop {\dot -} b, b \mathop {\dot -} a). 
\]
In other words, a subtraction is simply two rectified subtractions, and therefore by Theorem \ref{thm:composite}, has the same lower bound on the speed as either of the two rectified subtractions.

\subsection{Graphical representation for composite operations}
\label{sec:graph_rep}

\begin{figure}[H]
\begin{tabular}{c c c}
\begin{subfigure}[b]{0.3\textwidth}
\centering
\begin{circuitikz}
\draw 
(0,0) node[buffer port, fill=yellow] (id) {$\id$}
  (id.in) node[anchor=east] {$a$}
  (id.out) node[anchor=west] {$a$};
\end{circuitikz}
\caption{Identification.}
\end{subfigure}&
\begin{subfigure}[b]{0.3\textwidth}
\centering
\begin{circuitikz}
\draw 
(0,0) node[not port, fill=blue!40!] (inv) {$-1$}
  (inv.in) node[anchor=east] {$a$}
  (inv.out) node[anchor=west] {$1/a$};
   node{center};
\end{circuitikz}
\caption{Inversion.}
\end{subfigure}&
\begin{subfigure}[b]{0.3\textwidth}
\centering
\begin{circuitikz}
\draw 
(0,0) node[not port, fill=citrine] (mroot) {$m$}
  (mroot.in) node[anchor=east] {$a$}
  (mroot.out) node[anchor=west] {$\sqrt[m]{a}$};
\end{circuitikz}
\caption{$m$th root.}
\end{subfigure}\\
\begin{subfigure}[b]{0.3\textwidth}
\centering
\begin{circuitikz}
\draw 
(0,0) node[and port, fill=cyan!60!] (times) {$\times$}
  (times.in 1) node[anchor=east] {$a$}
  (times.in 2) node[anchor=east] {$b$}
  (times.out) node[anchor=west] {$a \times b$};
\end{circuitikz}
\caption{Multiplication.}
\end{subfigure}&
\begin{subfigure}[b]{0.3\textwidth}
\centering
\begin{circuitikz}
\draw 
(0,0) node[or port, fill=green] (plus) {$+$}
  (plus.in 1) node[anchor=east] {$a$}
  (plus.in 2) node[anchor=east] {$b$}
  (plus.out) node[anchor=west] {$a + b$};
\end{circuitikz}
\caption{Addition.}
\end{subfigure}&
\begin{subfigure}[b]{0.3\textwidth}
\centering
\begin{circuitikz} \draw
(0,0) node[xor port, fill=burntsienna] (abs) {$\abs{-}$}

(abs.in 1) node[anchor=east] {$a$}
(abs.in 2) node[anchor=east] {$b$}
(abs.out) node[anchor=west] {$\abs{a-b}$};

\end{circuitikz}
\caption{Absolute difference.}
\end{subfigure}\\
\begin{subfigure}[b]{0.3\textwidth}
\centering
\begin{circuitikz}
\draw 
(0,0) node[nor port, fill=bittersweet] (monus) {$\mathop {\dot -}$}
  (monus.in 1) node[anchor=east] {$a$}
  (monus.in 2) node[anchor=east] {$b$}
  (monus.out) node[anchor=west] {$a \mathop {\dot -} b$};
\end{circuitikz}
\caption{Rectified subtraction.}
\end{subfigure}&
\begin{subfigure}[b]{0.3\textwidth}
\centering
\begin{circuitikz} \draw
(0,0) node[op amp, noinv input up, fill=brightturquoise] (invr) {\circled{-1}}
(-2.5,0) node[circle, fill=yellow] (a) {$a$}
(-1.5,0.5) node[] (ap) {$a_p$}
(-1.5,-0.5) node[] (an) {$a_n$}
(a) -- (ap) -- (invr.+)
(a) -- (an) -- (invr.-)
(invr.out) node[right] {$\ds \begin{cases} \frac{1}{a_p}, &\mbox{if } a_p > 0, a_n=0, \\ 0, &\mbox{if } a_p = 0, a_n>0.\end{cases}$}
;
\end{circuitikz}
\caption{Partial real inversion.}
\end{subfigure}&
~
\end{tabular}
\caption{Elementary gates.}
\label{fig:elementary_gates}
\end{figure}
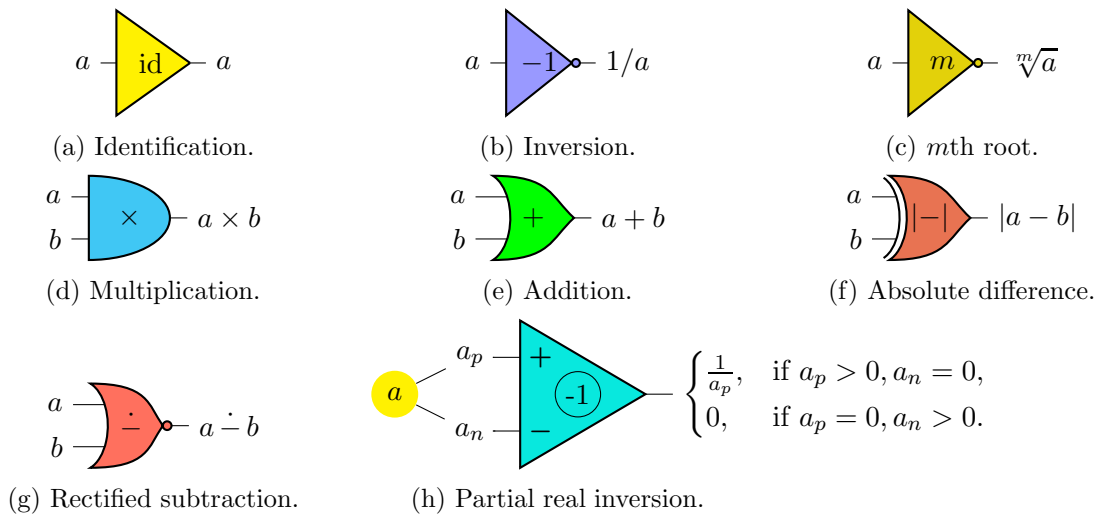

Combining elementary arithmetic computations into a composite computation has similarities with how logic gates are combined in a Boolean logic operation.  
Following this analogy, we can set up a visual aid where the elementary operations of section \ref{sec:elem_comp} are thought of ``elementary gates'' which can be combined in various different ways. 
We first give a catalogue of the elementary gates, shown in Figure \ref{fig:elementary_gates}.

\begin{figure}[H]
\centering
\begin{circuitikz} \draw
(3,1) node[nor port, fill=bittersweet] (monusa) {$\mathop {\dot -}$}
(3,-1) node[nor port, fill=bittersweet] (monusb) {$\mathop {\dot -}$}

(0,1.2) node[circle, fill=yellow] (a) {$a$}
(0,-1.2) node[circle, fill=blue!30!] (b) {$b$}

(monusa.out) node[anchor=west] {$a \mathop {\dot -} b$}
(monusb.out) node[anchor=west] {$b \mathop {\dot -} a$}

(a) -- (monusa.in 1)
(b) -- (monusa.in 2)
(a) -- (monusb.in 2)
(b) -- (monusb.in 1)
;
\end{circuitikz}
\caption{Subtraction.}
\label{fig:subtraction}
\end{figure}
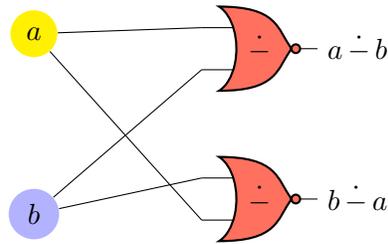

We can now put the gates together in different ways to generate the composite computations discussed in section \ref{sec:main_theory}. 
We start with division, shown in Figure \ref{fig:division} which only requires the elementary operations of identification, reciprocal, and multiplication. 
Figure \ref{fig:maximum} depicts a diagram of how the operation of taking the maximum of two non-negative reals can be constructed. 
Figure \ref{fig:subtraction} shows a pictorial way of doing subtraction by combining two rectified subtractions. 

We continue to use the diagrammatic representations of arithmetic computations in section \ref{sec:reals}. 
We show that each operation over reals can be decomposed into operations whose inputs and outputs are only non-negative reals. We give the diagrammatic decomposition into elementary arithmetic gates along with each description.

\begin{figure}[H]
\centering
\begin{circuitikz} \draw

(0.5,2) node[circle, fill=yellow] (a) {$a$}
(-1.2,0) node[circle, fill=green!40!] (b) {$b$}

(0,0) node[not port, fill=blue!40!] (inv) {$-1$}
(2,1) node[and port, fill=cyan!60!] (times) {$\times$}

(b) -- (inv.in)
(times.out) node[anchor=west] {$\ds \frac{a}{b}$}

(a) -- (times.in 1)
(inv.out) -- (times.in 2);
\end{circuitikz}
\caption{A graphical representation of division.}
\label{fig:division}
\end{figure}
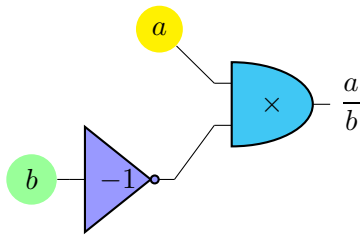

\begin{figure}[H]
\centering
\begin{circuitikz} \draw
(0,0) node[xor port, fill=burntsienna] (abs) {$\abs{-}$}
(2,0) node[or port, fill=green] (plusa) {$+$}
(4,0) node[or port, fill=green] (plusb) {$+$}
(6,0) node[and port, fill=cyan!60!] (times) {$\times$}

(-0.2,2) node[circle, fill=yellow] (a) {$a$}
(-0.2,-2) node[circle, fill=blue!30!] (b) {$b$}
(3.2,-1.5) node[circle, fill=brown!40!] (half) {$\frac12$}

(times.out) node[anchor=west] {$\ds \max\{a,b\}$}

(a) -- (abs.in 1)
(b) -- (abs.in 2)
(a) -- (plusa.in 1)
(b) -- (plusb.in 2)
(abs.out)  -- (plusa.in 2)
(plusa.out) -- (plusb.in 1)
(plusb.out) -- (times.in 1)

(half) -- (times.in 2) 
;
\end{circuitikz}
\caption{Maximum of two non-negative reals.}
\label{fig:maximum}
\end{figure}
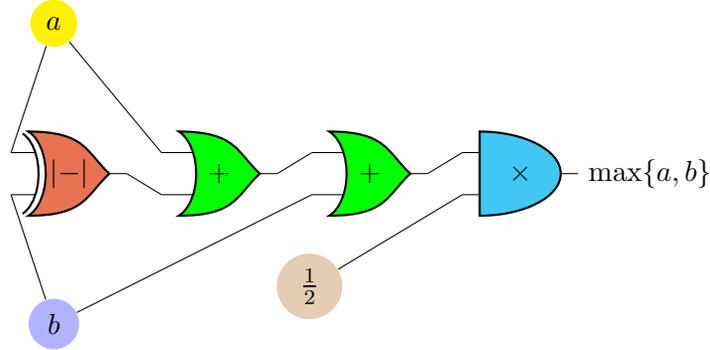

\section{Arithmetic operations on arbitrary real numbers}
\label{sec:reals}

The previous section discussed arithmetic operations with non-negative inputs. 
In the final construction on subtraction, the output is any real number--positive, negative, or zero--even though the inputs are positive. 
We needed `dual rail' to represent an arbitrary real number as a vector with two non-negative components. 
In this section, we continue to use the dual rail construction to represent an arbitrary real number which can be both an input and an output.

\subsubsection*{Addition and normalization}

We define addition in the usual way, for any $a, b \in \R$, 
\begin{align*}
a + b = (a_p,a_n) + (b_p, b_n) = (a_p + b_p, a_n + b_n). 
\end{align*}
The operation can be performed with ODEs/reaction networks by using two successive regular additions. 
Since the output of addition is not necessarily canonical, we define a normalization operation which brings the output into this form.

\[
(a_p, a_n) = 
\begin{cases}
(a_p - a_n,0)  & \mbox{if } a_p > a_n \\ 
(0,a_n - a_p) & \mbox{if } a_p \le a_n. 
\end{cases}
\]
The operation is performed using one dual rail subtraction. 

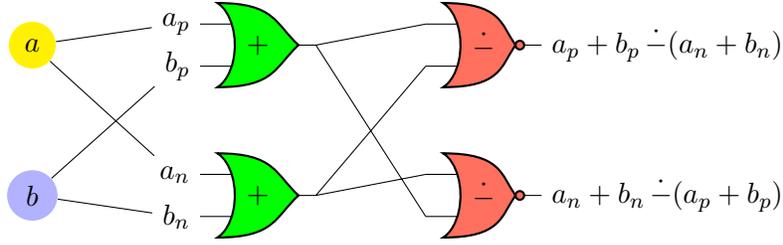
\begin{figure}[H]
\centering
\begin{circuitikz} \draw
(0,1) node[or port, fill=green] (plusp) {$+$}
  (plusp.in 1) node[anchor=east] (ghostap) {$a_p$}
  (plusp.in 2) node[anchor=east] (ghostbp) {$b_p$}
(0,-1) node[or port, fill=green] (plusn) {$+$}
  (plusn.in 1) node[anchor=east] (ghostan) {$a_n$}
  (plusn.in 2) node[anchor=east] (ghostbn) {$b_n$}
(3,1) node[nor port, fill=bittersweet] (monusa) {$\mathop {\dot -}$}
(3,-1) node[nor port, fill=bittersweet] (monusb) {$\mathop {\dot -}$}

(-3,1) node[circle, fill=yellow] (a) {$a$}
(-3,-1) node[circle, fill=blue!30!] (b) {$b$}

(monusa.out) node[anchor=west] {$a_p + b_p \mathop {\dot -} (a_n + b_n)$}
(monusb.out) node[anchor=west] {$a_n + b_n \mathop {\dot -} (a_p + b_p)$}

(a) -- (ghostap)
(a) -- (ghostan)
(b) -- (ghostbp)
(b) -- (ghostbn)
(plusp.out) -- (monusa.in 1)
(plusn.out) -- (monusa.in 2)
(plusp.out) -- (monusb.in 2)
(plusn.out) -- (monusb.in 1)
;
\end{circuitikz}
\caption{Real addition.}
\label{fig:realaddition}
\end{figure}

\subsubsection*{Additive inverse and subtraction}

The negative of the real number $a$ is defined via
\[
x = -a = -(a_p, a_n) = (a_n, a_p) = (x_p, x_n). 
\]
This can be performed using two assignments. 
The system of ODEs is the following:

\begin{align*}
\dot x_p &= a_n - x_p, \\
\dot x_n &= a_p - x_n. 
\end{align*}

Subtraction is defined as follows:
\begin{align*}
a - b = (a_p,a_n) - (b_p, b_n) = (a_p,a_n) + (- (b_p, b_n)) = (a_p,a_n) +  (b_n, b_p) = (a_p + b_n, a_n + b_p). 
\end{align*}

Thus subtraction can be performed using two regular additions followed by a normalization to bring the output into canonical form. 

\begin{figure}[H]
\centering
\begin{circuitikz} \draw
(3,1) node[buffer port, fill=yellow] (id1) {$\id$}
(3,-1) node[buffer port, fill=yellow] (id2) {$\id$}

(0,0) node[circle, fill=blue!30!] (a) {$a$}
(1,1) node[] (ap) {$a_p$}
(1,-1) node[] (an) {$a_n$}

(id1.out) node[anchor=west] {$a_n$}
(id2.out) node[anchor=west] {$a_p$}

(a) -- (ap) -- (id2.in)
(a) -- (an) -- (id1.in)
;
\end{circuitikz}
\caption{Real additive inverse.}
\label{fig:negative}
\end{figure}
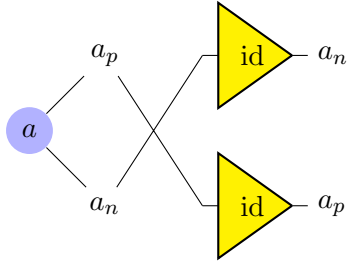

\subsubsection*{Multiplication}

Multiplication of arbitrary real numbers involves four multiplications and two additions of non-negative real numbers as shown by the following calculation: 
\begin{align*}
a \cdot b &= (a_p,a_n) \cdot (b_p, b_n) = (a_p - a_n) \cdot (b_p - b_n) \\
& = a_p b_p + a_n b_n - a_p b_n - a_n b_p = (a_p b_p + a_n b_n, a_p b_n + a_n b_p). 
\end{align*}
The normalization step is unnecessary since the output of multiplication is in canonical form whenever the inputs are in canonical form. Indeed, we have that 
\[
(a_p b_p + a_n b_n)(a_p b_n + a_n b_p) = (a_p^2+a_n^2)b_pb_n + a_pa_n(b_p^2+b_n^2) = 0, 
\]
where the last equality follows from $a_pa_n=0$ and $b_pb_n=0$. 

\begin{figure}[H]
\centering
\begin{circuitikz} \draw
(1,3) node[and port, fill=cyan!60!] (times1) {$\times$}
(1,1) node[and port, fill=cyan!60!] (times2) {$\times$}
(1,-1) node[and port, fill=cyan!60!] (times3) {$\times$}
(1,-3) node[and port, fill=cyan!60!] (times4) {$\times$}
(3,2) node[or port, fill=green] (plus1) {$+$}
(3,-2) node[or port, fill=green] (plus2) {$+$}

(-3,2) node[circle, fill=yellow] (a) {$a$}
(-3,-2) node[circle, fill=blue!30!] (b) {$b$}
(-2,3) node[] (ap) {$a_p$}
(-2,1) node[] (an) {$a_n$}
(-2,-1) node[] (bp) {$b_p$}
(-2,-3) node[] (bn) {$b_n$}

(plus1.out) node[anchor=west] {$a_p b_p + a_n  b_n$}
(plus2.out) node[anchor=west] {$a_p b_n + a_n  b_p$}

(a) -- (ap) -- (times1.in 1)
(ap) -- (times3.in 1)
(a) -- (an) -- (times2.in 1)
(an) -- (times4.in 1)
(b) -- (bp) -- (times1.in 2)
(bp) -- (times4.in 2)
(b) -- (bn) -- (times2.in 2)
(bn) -- (times3.in 2)
(times1.out) -- (plus1.in 1)
(times2.out) -- (plus1.in 2)
(times3.out) -- (plus2.in 1)
(times4.out) -- (plus2.in 2)
;
\end{circuitikz}
\caption{Real multiplication.}
\label{fig:realmultiplication}
\end{figure}
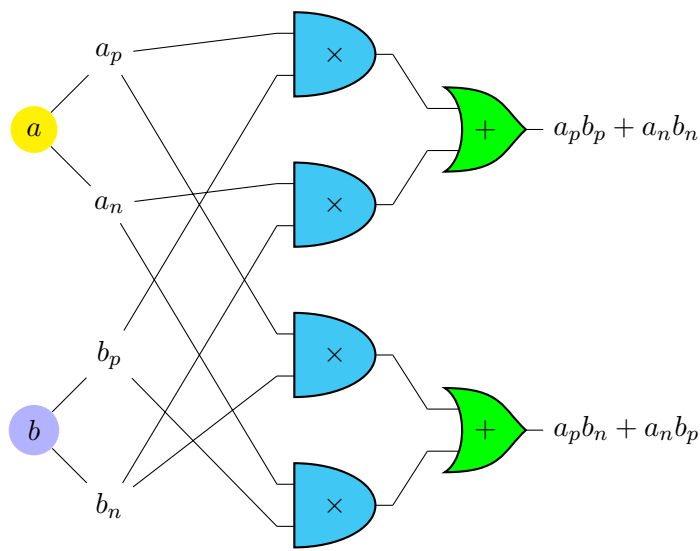

\subsubsection*{Inversion and Division}

Consider the nonzero input $0 \ne a = (a_p, a_n)$ in canonical form, i.e. exactly one of $a_p$ or $a_n$ is positive, so that $a_p + a_n > 0$ and $a_p a_n = 0$. 
The multiplicative inverse in canonical form is 
\[
a^{-1} = 
\begin{cases}
\left(\frac{1}{a_p},0\right)  & \mbox{if } a_p > 0, \\ 
\left(0, \frac{1}{a_n}\right) & \mbox{if } a_n > 0. 
\end{cases}
\]

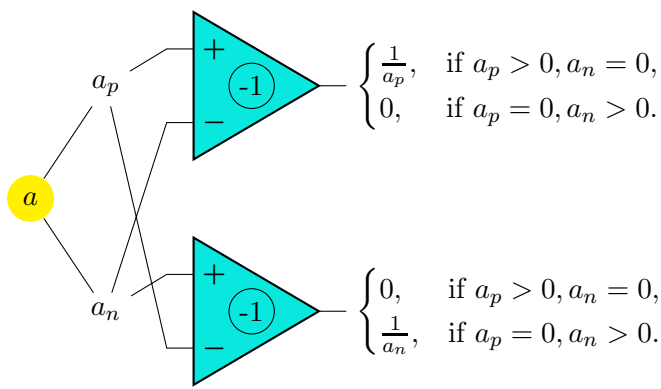
\begin{figure}[H]
\centering
\begin{circuitikz} \draw
(0,1.5) node[op amp, noinv input up, fill=brightturquoise] (invr1) {\circled{-1}}
(0,-1.5) node[op amp, noinv input up, fill=brightturquoise] (invr2) {\circled{-1}}

(-3,0) node[circle, fill=yellow] (a) {$a$}
(-2,1.5) node[] (ap) {$a_p$}
(-2,-1.5) node[] (an) {$a_n$}

(invr1.out) node[right] {$\ds \begin{cases} \frac{1}{a_p}, &\mbox{if } a_p > 0, a_n=0, \\ 0, &\mbox{if } a_p = 0, a_n>0.\end{cases}$}
(invr2.out) node[right] {$\ds \begin{cases} 0, &\mbox{if } a_p > 0, a_n=0, \\ \frac{1}{a_n}, &\mbox{if } a_p = 0, a_n>0.\end{cases}$}

(a) -- (ap) -- (invr1.+)
(ap) -- (invr2.-)
(a) -- (an) -- (invr1.-)
(an) -- (invr2.+)
;
\end{circuitikz}
\caption{Real inversion.}
\label{fig:realinverse}
\end{figure}

Finally, division $b/a$ is simply inversion of $a$ followed by multiplication of $b$ and $1/a$.

\section{Discussion}
\label{sec:discussion}

In this paper, we gave novel, explicit constructions for computing elementary operations at input-independent speed and showed how the elementary operations can be combined into composite computations also computed at input-independent speed. 
In particular, in the novel  constructions that we give, such as \eqref{net:inv} used for computing the reciprocal, we achieved a {\em decoupling} of concentrations and the speed of computing.  
This is in contrast to the naive constructions, such as \eqref{net:inv_naive} for the reciprocal.  
There, even the number of steps of the overall computation can affect the speed at which the individual steps are computed 
and therefore affect the speed of the overall computation.
This is because there could be complicated interactions between the number of steps in a given computation, the well-mixed property required for mass-action kinetics, the natural discreteness of the molecules, the volume of the vessel in which the reactions are taking place, and, therefore, the range of concentrations of the various species.  Such issues that arise from thinking about the physical implementation of chemical computation will be studied carefully in future work.

We assumed that underlying our algorithmic considerations is an idealized analog computer that is built from a perfectly well-behaved mass action system. 
We anticipate introducing realistic constraints in future research that will strengthen and build on the idealized setup developed here. 
We now discuss some future research directions that are natural offshoots of the present paper. 

Chemical reactions involving at most two molecules as reactants, sometimes referred to as bimolecular reactions, are considered the simplest to design and found most commonly in nature. The constructions given in this paper did not have this constraint and, in fact, did require multiple different reactant types with high molecularities. Is it possible to find constructions that compute arithmetic at input-independent speeds using only bimolecular reactants? 

In some of our constructions, an intermediate variable was allowed to increase to indefinitely high values. 
For instance when computing rectified subtraction \ref{eq:rectifiedsubtraction}, the concentration of the species $Y$ goes to infinity when $a^* = b^*$. 
This possibility of a variable increasing in an unbounded manner may be deemed unrealistic to implement using synthetic biochemistry, or at the very least is an energy-hungry construction. 
Is it possible to give an alternative construction for subtraction/rectified subtraction which works equally well algorithmically, has all the desired features, particularly including input-independent speed of computation, but also does not require any variable to increase in an unbounded manner? 

We assumed throughout that all reactions have the same rate constant. This poses significant, perhaps unrealistic, design challenges.
In the next improved iteration of construction, we may want to assume that the reactions occur at known, but not necessarily equal or controllable, rates. Given this significant constraint, how should the reaction networks for elementary operations and composite computations be modified to give an accurate output at input-independent speed? 

A bottleneck in the speed of computation for our constructions is a root of zero. Can this bottleneck be removed? Are there constructions for which the speed of computation is not lowered even when an elementary step requires computing a root of zero. 
Or is the jump in speed when computing a root of zero as opposed to a root of a positive number, not merely a bug but a positive feature? 
Perhaps we may be able to put to good use the discrete jump in speed when the input to an $m$-th root is either a zero or a positive number. 
To be concrete, consider the elementary operation of computing a square root. When the input is $\varepsilon > 0$, the speed of computation is at least $1$ while when the input is $0$, the speed is halved. 
Can this phenomenon be a basis for designing a ``zero detector'' or an ``equality detector''? 

Arithmetic only requires composing finitely many elementary operations. 
Computation of general functions may need the use of power series and other techniques that involve limits and convergence. 
Approximating the value of a general function, using a power series for example, may require computing an arbitrarily large, even if finite, number of elementary operations. 
To compute a power series, we need to consider the initial rate of convergence of each elementary operation, not merely its asymptotic rate. 
The `overhead' time required in each elementary operation may add up to a significant amount when considering a large, possibly infinite number of operations. 
Can our analysis be extended and what additional features or constraints are required when computing general functions? 

An alternative approach, that does not require using power series, may be to develop specific algorithms for computing commonly used functions such as exponentials and logarithms, based on mass-action kinetics and whose speeds are input-independent. 
A fruitful research direction may be to explore the computational reach of reaction network-based hardware by developing versatile algorithms that broaden the suite of possible computations at input-independent speeds. 
Conversely, hardware may impose additional constraints on the software thereby limiting computational accuracy or speed. 
Are there specific design aspects from the vantage of synthetic biochemistry or other practical engineering constraints that speed up or slow down computation? 
Multidisciplinary teams of researchers, representing synthetic biology, computer science as well as mathematics may be most poised to  answer related questions and make rapid progress. 

A research program could be founded on a broader territory of inquiry. 
As an instance, consider Boolean algebra, which is foundational to logical computations used in modern digital computers. One might ask: do there exist good, i.e. accurate and efficient, algorithms for performing Boolean algebra using reaction networks?

\section*{Acknowledgement}
We would like to thank Josef Hofbauer (University of Vienna) for suggesting the change of variables used in the proof of Lemma \ref{thm:system2_analysis}, which simplified the analysis. This work was supported by a grant from the Simons Foundation (905083, Anderson).  We also gratefully acknowledge NSF grant DMS-2051498.

\bibliographystyle{unsrt}
\bibliography{ACRN}

\appendix

\end{document}